\documentclass[12pt,leqno]{amsart}

\usepackage{amscd}
\usepackage{amsmath}
\usepackage{amsfonts}
\usepackage{amssymb}
\usepackage{url}

\title{Noether's problem for abelian extensions of cyclic $p$-groups}
\author{Ivo M. Michailov}
\address{Faculty of Mathematics and Informatics, Shumen University "Episkop Konstantin Preslavski", Universitetska str. 115, 9700 Shumen, Bulgaria}
\email{ivo\_michailov@yahoo.com}
\date{\today}
\thanks{This work is partially supported by a project No RD-05-275/15.01.2013 of Shumen University}
\begin{document}
\baselineskip 20pt
\begin{abstract}
Let $K$ be a field and $G$ be a finite group. Let $G$ act on the
rational function field $K(x(g):g\in G)$ by $K$ automorphisms
defined by $g\cdot x(h)=x(gh)$ for any $g,h\in G$. Denote by $K(G)$
the fixed field $K(x(g):g\in G)^G$. Noether's problem then asks
whether $K(G)$ is rational (i.e., purely transcendental) over $K$.
The first main result of this article is that $K(G)$ is rational
over $K$ for a certain class of $p$-groups having an abelian subgoup
of index $p$. The second main result is that $K(G)$ is rational over
$K$ for any group of order $p^5$ or $p^6$ ($p$ is an odd prime)
having an abelian normal subgroup such that its quotient group is
cyclic. (In both theorems we assume that if $char K\ne p$ then $K$
contains a primitive $p^e$-th root of unity, where $p^e$ is the
exponent of $G$.)
\end{abstract}

\maketitle
\newcommand{\Gal}{{\rm Gal}}
\newcommand{\Ker}{{\rm Ker}}
\newcommand{\GL}{{\rm GL}}
\newcommand{\Br}{{\rm Br}}
\newcommand{\lcm}{{\rm lcm}}
\newcommand{\ord}{{\rm ord}}
\renewcommand{\thefootnote}{\fnsymbol{footnote}}
\numberwithin{equation}{section}

\section{Introduction}
\label{1}

Let $K$ be any field. A field extension $L$ of $K$ is called
rational over $K$ (or $K$-rational, for short) if $L\simeq
K(x_1,\ldots,x_n)$ over $K$ for some integer $n$, with
$x_1,\ldots,x_n$ algebraically independent over $K$. Now let $G$ be
a finite group. Let $G$ act on the rational function field
$K(x(g):g\in G)$ by $K$ automorphisms defined by $g\cdot x(h)=x(gh)$
for any $g,h\in G$. Denote by $K(G)$ the fixed field $K(x(g):g\in
G)^G$. {\it Noether's problem} then asks whether $K(G)$ is rational
over $K$. This is related to the inverse Galois problem, to the
existence of generic $G$-Galois extensions over $K$, and to the
existence of versal $G$-torsors over $K$-rational field extensions
\cite[33.1, p.86]{Sw,Sa1,GMS}. Noether's problem for abelian groups
was studied extensively by Swan, Voskresenskii, Endo, Miyata and
Lenstra, etc. The reader is referred to Swan's paper for a survey of
this problem \cite{Sw}. Fischer's Theorem is a starting point of
investigating Noether's problem for finite abelian groups in
general.

\newtheorem{t1.1}{Theorem}[section]
\begin{t1.1}\label{t1.1}
{\rm (Fischer }\cite[Theorem 6.1]{Sw}{\rm )} Let $G$ be a finite
abelian group of exponent $e$. Assume that {\rm (i)} either char $K
= 0$ or char $K > 0$ with char $K\nmid e$, and {\rm (ii)} $K$
contains a primitive $e$-th root of unity. Then $K(G)$ is rational
over $K$.
\end{t1.1}

On the other hand, just a handful of results about Noether's problem
are obtained when the groups are non-abelian. This is the case even
when the group G is a $p$-group. The reader is referred to
\cite{CK,HuK,Ka1,Ka2,Ka3} for previous results of Noether's problem
for $p$-groups. The following theorem of Kang generalizes Fischer's
theorem for the metacyclic $p$-groups.

\newtheorem{t1.2}[t1.1]{Theorem}
\begin{t1.2}\label{t1.2}
{\rm (Kang}\cite[Theorem 1.5]{Ka1}{\rm )} Let $G$ be a metacyclic
$p$-group with exponent $p^e$, and let $K$ be any field such that
{\rm (i)} char $K = p$, or {\rm (ii)} char $K \ne p$ and $K$
contains a primitive $p^e$-th root of unity. Then $K(G)$ is rational
over $K$.
\end{t1.2}

The next job is to study Noether's problem for meta-abelian groups.
Three results due to Haeuslein, Hajja and Kang respectively are
known.

\newtheorem{t1.23}[t1.1]{Theorem}
\begin{t1.23}\label{t1.23}
{\rm (Haeuslein }\cite{Ha}{\rm )} Let $K$ be a field and $G$ be a
finite group. Assume that (i) $G$ contains an abelian normal
subgroup $H$ so that $G/H$ is cyclic of prime order $p$, (ii)
$\mathbb Z[\zeta_p]$ is a unique factorization domain, and (iii)
$\zeta_{p^e}\in K$ where $e$ is the exponent of $G$. If $G\to
\GL(V)$ is any finite-dimensional linear representation of $G$ over
$K$, then $K(V)^G$ is rational over $K$.
\end{t1.23}

\newtheorem{t1.24}[t1.1]{Theorem}
\begin{t1.24}\label{t1.24}
{\rm (Hajja }\cite{Haj}{\rm )} Let $K$ be a field and $G$ be a
finite group. Assume that (i) $G$ contains an abelian normal
subgroup $H$ so that $G/H$ is cyclic of order $n$, (ii) $\mathbb
Z[\zeta_n]$ is a unique factorization domain, and (iii) $K$ is
algebraically closed with $char K = 0$. If $G\to \GL(V)$ is any
finite-dimensional linear representation of $G$ over $K$, then
$K(V)^G$ is rational over $K$.
\end{t1.24}

\newtheorem{t1.22}[t1.1]{Theorem}
\begin{t1.22}\label{t1.22}
{\rm (}\cite[Theorem 1.4]{Ka3}{\rm )} Let $K$ be a field and $G$ be
a finite group. Assume that {\rm (i)} $G$ contains an abelian normal
subgroup $H$ so that $G/H$ is cyclic of order $n$, {\rm (ii)}
$\mathbb Z[\zeta_n]$ is a unique factorization domain, and {\rm
(iii)} $\zeta_{e}\in K$ where $e$ is the exponent of $G$. If
$G\rightarrow \GL(V)$ is any finite-dimensional linear
representation of $G$ over $K$, then $K(V)^G$ is rational over $K$.
\end{t1.22}

Note that those integers $n$ for which $\mathbb Z[\zeta_n]$ is a
unique factorization domain are determined by Masley and Montgomery.

\newtheorem{t1.25}[t1.1]{Theorem}
\begin{t1.25}\label{t1.25}
{\rm (Masley and Montgomery }\cite{MM}{\rm )}  $\mathbb Z[\zeta_n]$
is a unique factorization domain if and only if $1\leq n\leq 22$, or
$n = 24, 25, 26, 27, 28, 30, 32, 33, 34, 35, 36, 38, 40, 42, 45,
48,$ $50, 54, 60, 66, 70, 84, 90$.
\end{t1.25}

Therefore, Theorem \ref{t1.23} holds only for primes $p$ such that
$1\leq p\leq 19$. One of the goals of our paper is to show that the
this condition can be waived, under some additional assumptions
regarding the structure of the abelian subgroup $H$.

Consider the following situation. Let $G$ be a group of order $p^n$
for $n\geq 2$ with an abelian subgroup $H$ of order $p^{n-1}$.
Bender \cite{Be2} determined some interesting properties of these
groups. We study further the case when the $p$-th lower central
subgroup $G_{(p)}$ is trivial. (Recall that $G_{(0)}=G$ and
$G_{(i)}=[G,G_{(i-1)}]$ for $i\geq 1$ are called the lower central
series.) For our purposes we need to classify with generators and
relations these groups. We achieve this in the following lemma.

\newtheorem{lemma}[t1.1]{Lemma}
\begin{lemma}\label{lemma}
Let $G$ be a group of order $p^n$ for $n\geq 2$ with an abelian
subgroup $H$ of order $p^{n-1}$. Choose any $\alpha\in G$ such that
$\alpha$ generates $G/H$, i.e., $\alpha\notin H,\alpha^p\in H$.
Denote $H(p)=\{h\in H: h^p=1,h\notin H^p\}$, and assume that
$[H(p),\alpha]\subset H(p)$. Assume also that the $p$-th lower
central subgroup $G_{(p)}$ is trivial. (Recall that $G_{(0)}=G$ and
$G_{(i)}=[G,G_{(i-1)}]$ for $i\geq 1$ are called the lower central
series.) Then $H$ is a direct product of normal subgroups of $G$
that are of the following four types:
\begin{enumerate}
    \item $(C_p)^s$ for some $s\geq 1$. There exist generators $\alpha_1,\dots,\alpha_s$ of
$(C_p)^s$, such that $[\alpha_j,\alpha]=\alpha_{j+1}$ for $1\leq
j\leq s-1$ and $\alpha_s\in Z(G)$.
    \item $C_{p^a}$ for some $a\geq 1$. There exists a
    generator $\beta$ of $C_{p^a}$ such that $[\beta,\alpha]=\beta^{bp^{a-1}}$ for some $b:0\leq b\leq p-1$.
    \item $C_{p^{a_1}}\times
C_{p^{a_2}}\times\cdots\times C_{p^{a_k}}\times(C_p)^s$ for some
$k\geq 1,a_i\geq 2,s\geq 1$. There exist generators
$\alpha_{11},\alpha_{21},\dots,\alpha_{k1}$ of $C_{p^{a_1}}\times
C_{p^{a_2}}\times\cdots\times C_{p^{a_k}}$ such that
$[\alpha_{i,1},\alpha]=\alpha_{i+1,1}^{p^{a_{i+1}-1}}\in Z(G)$ for
$i=1,\dots,k-1$. There also exist
    generators
$\alpha_{k,2},\dots,\alpha_{k,s+1}$ of $(C_p)^s$, such that
$[\alpha_{k,j},\alpha]=\alpha_{k,j+1}$ for $1\leq j\leq s$ and
$\alpha_{k,s+1}\in Z(G)$.
    \item $C_{p^{a_1}}\times
C_{p^{a_2}}\times\cdots\times C_{p^{a_k}}$ for some $k\geq 2,a_i\geq
2$. For any $i:1\leq i\leq k$ there exists a generator
$\alpha_{i,1}$ of the factor $C_{p^{a_i}}$, such that
$[\alpha_{i,1},\alpha]=\alpha_{i+1,1}^{p^{a_i-1}}\in Z(G)$ and
$[\alpha_{k,1},\alpha]\in\langle\alpha_{1,1}^{p^{a_1-1}},\dots,\alpha_{k,1}^{p^{a_{k}-1}}\rangle$.
\end{enumerate}
\end{lemma}

The first main result of this paper is the following theorem which
generalizes Theorem \ref{t1.23}.

\newtheorem{t1.4}[t1.1]{Theorem}
\begin{t1.4}\label{t1.4}
Let $G$ be a group of order $p^n$ for $n\geq 2$ with an abelian
subgroup $H$ of order $p^{n-1}$, and let $G$ be of exponent $p^e$.
Choose any $\alpha\in G$ such that $\alpha$ generates $G/H$, i.e.,
$\alpha\notin H,\alpha^p\in H$. Denote $H(p)=\{h\in H: h^p=1,h\notin
H^p\}$, and assume that $[H(p),\alpha]\subset H(p)$. Denote by
$G_{(i)}=[G,G_{(i-1)}]$ the lower central series for $i\geq 1$ and
$G_{(0)}=G$. Let the $p$-th lower central subgroup $G_{(p)}$ be
trivial. Assume that {\rm (i)} char $K = p>0$, or {\rm (ii)} char $K
\ne p$ and $K$ contains a primitive $p^e$-th root of unity. Then
$K(G)$ is rational over $K$.
\end{t1.4}

The key idea to prove Theorem \ref{t1.4} is to find a faithful
$G$-subspace $W$ of the regular representation space
$\bigoplus_{g\in G} K\cdot x(g)$ and to show that $W^G$ is rational
over $K$. The subspace $W$ is obtained as an induced representation
from $H$ by applying Lemma \ref{lemma}. (A particular case we proved
in the preprint \cite{Mi}.)

The next goal of our article is to study Noether's problem for some
groups of orders $p^5$ and $p^6$ for any odd prime $p$. We use the
list of generators and relations for these groups, given by James
\cite{Ja}. It is known that $K(G)$ is always rational if $G$ is a
$p$-group of order $\leq p^4$ and $\zeta_{e}\in K$ where $e$ is the
exponent of $G$ (see \cite{CK}). However, in \cite{HoK} is shown
that there exists a group $G$ of order $p^5$ such that $\mathbb
C(G)$ is not rational over $\mathbb C$.

The second main result of this article is the following rationality
criterion for the groups of orders $p^5$ and $p^6$, having an
abelian normal subgroup such that its quotient group is cyclic.

\newtheorem{t1.5}[t1.1]{Theorem}
\begin{t1.5}\label{t1.5}
Let $G$ be a group of order $p^n$ for $n\leq 6$ with an abelian
normal subgroup $H$, such that $G/H$ is cyclic. Let $G$ be of
exponent $p^e$. Assume that {\rm (i)} char $K = p>0$, or {\rm (ii)}
char $K \ne p$ and $K$ contains a primitive $p^e$-th root of unity.
Then $K(G)$ is rational over $K$.
\end{t1.5}

We do not know whether Theorem \ref{t1.5} holds for any $n\geq 7$.
However, we should not ''over-generalize'' Theorem \ref{t1.5} to the
case of any meta-abelian group because of the following theorem of
Saltman.

\newtheorem{t1.3}[t1.1]{Theorem}
\begin{t1.3}\label{t1.3}
{\rm (Saltman }\cite{Sa2}{\rm )} For any prime number $p$ and for
any field $K$ with char $K \ne p$ (in particular, $K$ may be an
algebraically closed field), there is a meta-abelian $p$-group $G$
of order $p^9$ such that $K(G)$ is not rational over $K$.
\end{t1.3}

We organize this paper as follows. We recall some preliminaries in
Section \ref{2} that will be used in the proofs of Theorems
\ref{t1.4} and \ref{t1.5}. There we also prove Lemma \ref{l2.8}
which is a generalization of Kang's argument from \cite[Case 5, Step
II]{Ka2}. In Section \ref{3} we prove Lemma \ref{lemma} which is of
independent interest, since it provides a list of generators and
relations for any $p$-group $G$ having an abelian subgroup $H$ of
index $p$, provided that $[H(p),\alpha]\subset H(p)$ and
$G_{(p)}=1$. Our main results -- Theorems \ref{t1.4} and \ref{t1.5}
-- we prove in Sections \ref{4} and \ref{5} respectively.

\section{Preliminaries}
\label{2}

We list several results which will be used in the sequel.

\newtheorem{t2.1}{Theorem}[section]
\begin{t2.1}\label{t2.1}
{\rm (}\cite[Theorem 1]{HK}{\rm )} Let $G$ be a finite group acting
on $L(x_1,\dots,x_m)$, the rational function field of $m$ variables
over a field $L$ such that
\begin{description}
    \item [(i)] for any $\sigma\in G, \sigma(L)\subset L;$
    \item [(ii)] the restriction of the action of $G$ to $L$ is
    faithful;
    \item [(iii)] for any $\sigma\in G$,
    \begin{equation*}
\begin{pmatrix}
\sigma(x_1)\\
\vdots\\
\sigma(x_m)\\
\end{pmatrix}
=A(\sigma)\begin{pmatrix}
x_1\\
\vdots\\
x_m\\
\end{pmatrix}
+B(\sigma)
\end{equation*}
where $A(\sigma)\in\GL_m(L)$ and $B(\sigma)$ is $m\times 1$ matrix
over $L$. Then there exist $z_1,\dots,z_m\in L(x_1,\dots,x_m)$ so
that $L(x_1,\dots,x_m)^G=L^G(z_1,\dots,z_m)$ and $\sigma(z_i)=z_i$
for any $\sigma\in G$, any $1\leq i\leq m$.
\end{description}
\end{t2.1}

\newtheorem{t2.2}[t2.1]{Theorem}
\begin{t2.2}\label{t2.2}
{\rm (}\cite[Theorem 3.1]{AHK}{\rm )} Let $G$ be a finite group
acting on $L(x)$, the rational function field of one variable over a
field $L$. Assume that, for any $\sigma\in G,\sigma(L)\subset L$ and
$\sigma(x)=a_\sigma x+b_\sigma$ for any $a_\sigma,b_\sigma\in L$
with $a_\sigma\ne 0$. Then $L(x)^G=L^G(z)$ for some $z\in L[x]$.
\end{t2.2}

\newtheorem{t2.3}[t2.1]{Theorem}
\begin{t2.3}\label{t2.3}
{\rm (}\cite[Theorem 1.7]{CK}{\rm )} If $char K=p>0$ and $\widetilde
G$ is a finite $p$-group, then $K(G)$ is rational over $K$.
\end{t2.3}

The following Lemma can be extracted from some proofs in
\cite{Ka2,HuK}.

\newtheorem{l2.7}[t2.1]{Lemma}
\begin{l2.7}\label{l2.7}
Let $\langle\tau\rangle$ be a cyclic group of order $n>1$, acting on
$K(v_1,\dots,v_{n-1})$, the rational function field of $n-1$
variables over a field $K$ such that
\begin{eqnarray*}
\tau&:&v_1\mapsto v_2\mapsto\cdots\mapsto v_{n-1}\mapsto (v_1\cdots
v_{n-1})^{-1}\mapsto v_1.
\end{eqnarray*}
If $K$ contains a primitive $n$-th root of unity $\xi$, then
$K(v_1,\dots,v_{n-1})=K(s_1,\dots,s_{n-1})$ where $\tau:s_i\mapsto
\xi^is_i$ for $1\leq i\leq n-1$.
\end{l2.7}
\begin{proof}
Define $w_0=1+v_1+v_1v_2+\cdots+v_1v_2\cdots
v_{n-1},w_1=(1/w_0)-1/n,w_{i+1}=(v_1v_2\cdots v_i/w_0)-1/n$ for
$1\leq i\leq n-1$. Thus $K(v_1,\dots,v_{n-1})=K(w_1,\dots,w_n)$ with
$w_1+w_2+\cdots+w_n=0$ and
\begin{eqnarray*}
\tau&:&w_1\mapsto w_2\mapsto\cdots\mapsto w_{n-1}\mapsto w_n\mapsto
w_1.
\end{eqnarray*}
Define $s_i=\sum_{1\leq j\leq n}\xi^{-ij}w_j$ for $1\leq i\leq n-1$.
Then $K(w_1,\dots,w_n)=K(s_1,\dots,s_{n-1})$ and $\tau:s_i\mapsto
\xi^is_i$ for $1\leq i\leq n-1$.
\end{proof}

Moreover, we are now going to generalize Kang's argument from
\cite[Case 5, Step II]{Ka2}, obtaining the following Lemma which
plays an important role in our work.

\newtheorem{l2.8}[t2.1]{Lemma}
\begin{l2.8}\label{l2.8}
Let $k>1$, let $p$ be any prime and let $\langle\alpha\rangle$ be a
cyclic group of order $p$, acting on
$K(y_{1i},y_{2i}\dots,y_{ki}:1\leq i\leq p-1)$, the rational
function field of $k(p-1)$ variables over a field $K$ such that
\begin{align*} \alpha\ :\ &y_{j1}\mapsto
y_{j2}\mapsto\cdots\mapsto y_{jp-1}\mapsto (y_{j1}y_{j2}\cdots
y_{jp-1})^{-1},\ &\text{for}\ 1\leq j\leq k.
\end{align*}
Assume that $K(v_{1i},v_{2i}\dots,v_{ki}:1\leq i\leq
p-1)=K(y_{1i},y_{2i}\dots,y_{ki}:1\leq i\leq p-1)$ where for any
$j:1\leq j\leq k$ and for any $i:1\leq i\leq p-1$ the variable
$v_{ji}$ is a monomial in the variables $y_{1i},y_{2i}\dots,y_{ki}$.
Assume also that the action of $\alpha$ on
$K(v_{1i},v_{2i}\dots,v_{ki}:1\leq i\leq p-1)$ is given by
{\allowdisplaybreaks\begin{align*} \alpha\ :\ &v_{j1}\mapsto
v_{j1}v_{j2}^p,~ v_{j2}\mapsto v_{j3}\mapsto\cdots\mapsto
v_{jp-1}\mapsto
A_j\cdot(v_{j1}v_{j2}^{p-1}v_{j3}^{p-2}\cdots v_{jp-1}^2)^{-1},\\
\ &\text{for}\ 1\leq j\leq k,
\end{align*}}
where $A_j$ is some monomial in $v_{1i},\dots,v_{j-1i}$ for $2\leq
j\leq k$ and $A_1=1$. If $K$ contains a primitive $p$-th root of
unity $\zeta$, then $K(v_{1i},v_{2i}\dots,v_{ki}:1\leq i\leq
p-1)=K(s_{1i},s_{2i}\dots,s_{ki}:1\leq i\leq p-1)$ where
$\alpha:s_{ji}\mapsto \zeta^is_{ji}$ for $1\leq j\leq k,1\leq i\leq
p-1$.
\end{l2.8}
\begin{proof}
We write the additive version of the multiplication action of
$\alpha$, i.e., consider the $\mathbb Z[\pi]$-module
$M=\bigoplus_{1\leq m\leq k}(\oplus_{1\leq i\leq p-1}\mathbb Z\cdot
v_{mi})$, where $\pi=\langle\alpha\rangle$. Denote the submodules
$M_j=\bigoplus_{1\leq m\leq j}(\oplus_{1\leq i\leq p-1}\mathbb
Z\cdot v_{mi})$ for $1\leq j\leq k$. Thus $\alpha$ has the following
additive action {\allowdisplaybreaks\begin{align*}
\alpha\ :\ &v_{j1}\mapsto v_{j1}+pv_{j2},~\\
 &v_{j2}\mapsto
v_{j3}\mapsto\cdots\mapsto v_{jp-1}\mapsto
A_j-v_{j1}-(p-1)v_{j2}-(p-2)v_{j3}-\cdots -2v_{jp-1},
\end{align*}}
where $A_j\in M_{j-1}$.

By Lemma \ref{l2.7}, $M_1$ is isomorphic to the $\mathbb
Z[\pi]$-module $N=\oplus_{1\leq i\leq p-1}\mathbb Z\cdot u_i$ where
$u_1=v_{12},u_i=\alpha^{i-1}\cdot v_{12}$ for $2\leq i\leq p-1$, and
\begin{align*}
\alpha\ :\  &u_1\mapsto u_2\mapsto\cdots\mapsto u_{p-1}\mapsto
-u_1-u_2-\cdots-u_{p-1}\mapsto u_1.
\end{align*}

Let $\Phi_p(T)\in\mathbb Z[T]$ be the $p$-th cyclotomic polynomial.
Since $\mathbb Z[\pi]\simeq\mathbb Z[T]/(T^p-1)$, we find that
$\mathbb Z[\pi]/\Phi_p(\alpha)\simeq \mathbb Z[T]/\Phi_p(T)\simeq
\mathbb Z[\omega]$, the ring of $p$-th cyclotomic integer. As
$\Phi_p(\alpha)\cdot x=0$ for any $x\in N$, the $\mathbb
Z[\pi]$-module $N$ can be regarded as a $\mathbb Z[\omega]$-module
through the morphism $\mathbb Z[\pi]\to\mathbb
Z[\pi]/\Phi_p(\alpha)$. When $N$ is regarded as a $\mathbb
Z[\omega]$-module, $N\simeq\mathbb Z[\omega]$ the rank-one free
$\mathbb Z[\omega]$-module.

We claim that $M$ itself can be regarded as a $\mathbb
Z[\omega]$-module, i.e., $\Phi_p(\alpha)\cdot M=0$.

Return to the multiplicative notations. Note that all $v_{ji}$'s are
monomials in $y_{ji}$'s. The action of $\alpha$ on $y_{ji}$ given in
the statement satisfies the relation $\prod_{0\leq m\leq
p-1}\alpha^m(y_{ji})=1$ for any $1\leq j\leq k,1\leq i\leq p-1$.
Using the additive notations, we get $\Phi_p(\alpha)\cdot y_{ji}=0$.
Hence $\Phi_p(\alpha)\cdot M=0$.

Define $M'=M/M_{k-1}$. It follows that we have a short exact
sequence of $\mathbb Z[\pi]$-modules
\begin{equation}\label{e2.3}
0\to M_{k-1}\to M\to M'\to 0.
\end{equation}
Since $M$ is a $\mathbb Z[\omega]$-module, \eqref{e2.3} is a short
exact sequence of $\mathbb Z[\omega]$-modules. Proceeding by
induction, we obtain that $M$ is a direct sum of free $\mathbb
Z[\omega]$-modules isomorphic to $N$. Therefore,
$M\simeq\oplus_{1\leq j\leq k}N_j$, where $N_j\simeq N$ is a free
$\mathbb Z[\omega]$-module, and so a $\mathbb Z[\pi]$-module also
(for $1\leq j\leq k$).

Finally, we interpret the additive version of $M\simeq\oplus_{1\leq
j\leq k}N_j\simeq N^k$ it terms of the multiplicative version as
follows: There exist $w_{ji}$ that are monomials in $v_{ji}$ for
$1\leq j\leq k,1\leq i\leq p-1$ such that $K(w_{ji})=K(v_{ji})$ and
$\alpha$ acts as
\begin{align*}
\alpha\ :\ &w_{j1}\mapsto w_{j2}\mapsto\cdots\mapsto w_{jp-1}\mapsto
(w_{j1}w_{j2}\dots w_{jp-1})^{-1}\ \text{for}\ 1\leq j\leq k.
\end{align*}
According to Lemma \ref{l2.7}, the above action can be linearized as
pointed in the statement.
\end{proof}

Now, let $G$ be any metacyclic $p$-group generated by two elements
$\sigma$ and $\tau$ with relations
$\sigma^{p^a}=1,\tau^{p^b}=\sigma^{p^c}$ and
$\tau^{-1}\sigma\tau=\sigma^{\varepsilon+\delta p^r}$ where
$\varepsilon=1$ if $p$ is odd, $\varepsilon=\pm 1$ if $p=2$,
$\delta=0,1$ and $a,b,c,r\geq 0$ are subject to some restrictions.
For the the description of these restrictions see e.g. \cite[p.
564]{Ka1}.

\newtheorem{t2.6}[t2.1]{Theorem}
\begin{t2.6}\label{t2.6}
{\rm (Kang }\cite[Theorem 4.1]{Ka1}{\rm )} Let $p$ be a prime
number, $m,n$ and $r$ are positive integers, $k=1+p^r$ if $(p,r)\ne
(2,1)$ (resp. $k=-1+2^r$ with $r\geq 2$). Let $G$ be a split
metacyclic $p$-group of order $p^{m+n}$ and exponent $p^e$ defined
by $G=\langle\sigma,\tau:
\sigma^{p^m}=\tau^{p^n}=1,\tau^{-1}\sigma\tau=\sigma^k\rangle$. Let
$K$ be any field such that $char K\ne p$ and $K$ contains a
primitive $p^e$-th root of unity, and let $\zeta$ be a primitive
$p^m$-th root of unity. Then $K(x_0,x_1,\dots,x_{p^n-1})^G$ is
rational over $K$, where $G$ acts on $x_0,\dots,x_{p^n-1}$ by
\begin{eqnarray*}
\sigma&:&x_i\mapsto \zeta^{k^i}x_i,\\
\tau&:&x_0\mapsto x_1\mapsto\cdots\mapsto x_{p^n-1}\mapsto x_0.
\end{eqnarray*}
\end{t2.6}

\section{Proof of Lemma \ref{lemma}}
\label{3}

It is well known that $H$ is a normal subgroup of $G$. We divide the
proof into several steps.

\emph{Step I.} Let $\beta_1$ be any element of $H$ that is not
central. Since $G_{(p)}=\{1\}$, there exist
$\beta_2,\dots,\beta_k\in H$ for some $k:2\leq k\leq p$ such that
$[\beta_j,\alpha]=\beta_{j+1}$, where $1\leq j\leq k-1$ and
$\beta_k\ne 1$ is central. We are going to show now that the order
of $\beta_2$ is not greater than $p$. In particular, from the
multiplication rule $[a,\alpha][b,\alpha]=[ab,\alpha]$ (for any
$a,b\in H$) it follows that all $p$-th powers are contained in the
center of $G$.

From $[\beta_j,\alpha]=\beta_{j+1}$ it follows the well known
formula
\begin{equation}\label{e3.1}
\alpha^{-p}\beta_1\alpha^p=\beta_1\beta_2^{\binom{p}{1}}\beta_3^{\binom{p}{2}}\cdots
\beta_p^{\binom{p}{p-1}}\beta_{p+1},
\end{equation}
where we put $\beta_{k+1}=\cdots=\beta_{p+1}=1$. Since $\alpha^p$ is
in $H$, we obtain the formula
$$\beta_2^{\binom{p}{1}}\beta_3^{\binom{p}{2}}\cdots
\beta_k^{\binom{p}{k-1}}=1.$$ Hence $(\beta_2\cdot\prod_{j\ne
2}\beta_j^{a_j})^p=1$ for some integers $a_j$. It is not hard to see
that this identity is impossible if the order of $\beta_2$ is
greater than $p$. Indeed, if $\ell=\max\{j:\beta_j^p\ne 1\}$, then
$\beta_\ell^p$ is in the subgroup generated by
$\beta_2^p,\dots,\beta_{\ell-1}^p$. Therefore
$[\beta_\ell^p,\alpha]=[\beta_2^{b_2p}\cdots\beta_{\ell-1}^{b_{\ell-1}p},\alpha]=\beta_3^{b_2p}\cdots\beta_\ell^{b_{\ell-1}p}\ne
1$ for some $b_2,\dots,b_{\ell-1}\in\mathbb Z_p$. On the other hand,
$[\beta_\ell^p,\alpha]=\beta_{\ell+1}^p=1$, which is a
contradiction.

\emph{Step II.} Let us write the decomposition of $H$ as a direct
product of cyclic subgroups (not necessarily normal in $G$):
$H\simeq (C_p)^t\times C_{p^{a_1}}\times C_{p^{a_2}}\times\cdots
C_{p^{a_s}}$ for $0\leq t, 2\leq a_1\leq a_2\leq\cdots\leq a_s$.
Choose a generator $\alpha_{11}\in C_{p^{a_1}}$. Since
$G_{(p)}=\{1\}$, there exist $\alpha_{12},\dots,\alpha_{1k}\in H$
for some $k:2\leq k\leq p$ such that
$[\alpha_{1j},\alpha]=\alpha_{1j+1}$, where $1\leq j\leq k-1$ and
$\alpha_{1k}\ne 1$ is central. From Step I it follows that the order
of $\alpha_{12}$ is not greater than $p$. We are going to define a
normal subgroup of $G$ which depends on the nature of the element
$\alpha_{12}$. We will denote it by
$\langle\langle\alpha_{11}\rangle\rangle$, and call it \emph{the
commutator chain of} $\alpha_{11}$. Simultaneously, we will define a
complement in $H$ denoted by
$\overline{\langle\langle\alpha_{11}\rangle\rangle}$.

\emph{Case II.1.} Let $\alpha_{12}=\alpha_{11}^{p^{a_1-1}c_1}$ for
some $c_1:0\leq c_1\leq p-1$. Define
$\langle\langle\alpha_{11}\rangle\rangle=\langle\alpha_{11}\rangle$.
Clearly, $\langle\langle\alpha_{11}\rangle\rangle$ is a normal
subgroup of type (2). Define
$\overline{\langle\langle\alpha_{11}\rangle\rangle}=(C_p)^t\cdot
\langle\alpha_{21},\dots,\alpha_{s1}\rangle$.

\emph{Case II.2.} Let $\alpha_{12}\notin H^p$. According to the
assumption in the statement of our Lemma, $[H(p),\alpha]\cap
H^p=\{1\}$, we have $\alpha_{1j}\notin H^p$ for all $j$. Define
$\langle\langle\alpha_{11}\rangle\rangle=\langle\alpha_{11},\dots,\alpha_{1k}\rangle$.
Therefore $\langle\langle\alpha_{11}\rangle\rangle\simeq
C_{p^{a_1}}\times (C_p)^{k-1}$ is a normal subgroup of type (3).
Define
$\overline{\langle\langle\alpha_{11}\rangle\rangle}=(C_p)^{t-k+1}\cdot
\langle\alpha_{21},\dots,\alpha_{s1}\rangle$, where $(C_p)^{t-k+1}$
is the complement of $(C_p)^{k-1}$ in $(C_p)^t$.

\emph{Case II.3.} Let $\alpha_{12}\in H^p$. Then
$\alpha_{12}=\prod_{i\in A}\alpha_{i1}^{p^{a_i-1}d_i}$, where
$A\subset\{1,2,\dots,s\},1\leq d_i\leq p-1$. Put $i_0=\min\{i\in
A\}$.

If $i_0=1$, then $\alpha_{12}=(\alpha_{11}^{d_1}\prod_{i\in A,i\ne
1}\alpha_{i1}^{p^{a_i-a_1}d_i})^{p^{a_1-1}}$. Now, we can replace
the generator $\alpha_{11}$ with
$\alpha_{11}'=\alpha_{11}^{d_1}\prod_{i\in A,i\ne
1}\alpha_{i1}^{p^{a_i-a_1}d_i}$. Clearly,
$\ord(\alpha_{11}')=\ord(\alpha_{11})$ and
$[\alpha_{11}',\alpha]\in\langle\alpha_{11}'\rangle$, so this case
is reduced to Case I.

If $i_0>1$, then $\alpha_{12}=(\alpha_{i_01}^{d_{i_0}}\prod_{i\in
A,i\ne i_0}\alpha_{i1}^{p^{a_i-a_{i_0}}d_i})^{p^{a_{i_0}-1}}$. We
can replace the generator $\alpha_{i_01}$ with
$\alpha_{i_01}'=\alpha_{i_01}^{d_{i_0}}\prod_{i\in A,i\ne
i_0}\alpha_{i1}^{p^{a_i-a_{i_0}}d_i}$. Clearly,
$\ord(\alpha_{i_01}')=\ord(\alpha_{i_01})$ and
$\alpha_{i_01}'^{p^{a_{i_0}-1}}=\alpha_{12}$.

For abuse of notation we will assume henceforth that $i_0=2$ and
$\alpha_{21}^{p^{a_2-1}}=\alpha_{12}$. Consider
$\alpha_{22}=[\alpha_{21},\alpha]$. We have three possibilities now.

\emph{Subcase II.3.1.}
$\alpha_{22}\in\langle\alpha_{11}^{p^{a_1-1}},\alpha_{21}^{p^{a_1-1}}\rangle$.
Define
$\langle\langle\alpha_{11}\rangle\rangle=\langle\alpha_{11},\alpha_{21}\rangle$.
Therefore $\langle\langle\alpha_{11}\rangle\rangle\simeq
C_{p^{a_1}}\times C_{p^{a_2}}$ is a normal subgroup of type (4).

\emph{Subcase II.3.2.} $\alpha_{22}\notin H^p$. Then there exist
$\alpha_{22},\dots,\alpha_{2\ell}\in H$ for some $\ell:2\leq
\ell\leq p$ such that $[\alpha_{2j},\alpha]=\alpha_{2j+1}$, where
$1\leq j\leq \ell-1$ and $\alpha_{2\ell}\ne 1$ is central. Define
$\langle\langle\alpha_{11}\rangle\rangle=\langle\alpha_{11},\alpha_{21},\alpha_{22},\dots,\alpha_{2\ell}\rangle$.
Therefore $\langle\langle\alpha_{11}\rangle\rangle\simeq
C_{p^{a_1}}\times C_{p^{a_2}}\times (C_p)^{\ell-1}$ is a normal
subgroup of type (3).

\emph{Subcase II.3.3.} $\alpha_{22}\in H^p$. According to the
observations we have just made, this subcase leads to the following
two final possibilities.

\emph{Sub-subcase II.3.3.1.}
$\alpha_{22}=\alpha_{31}^{p^{a_3-1}},\dots,\alpha_{r-12}=\alpha_{r1}^{p^{a_r-1}},\alpha_{r2}\in
\langle\alpha_{11}^{p^{a_1-1}},\dots,\alpha_{r1}^{p^{a_r-1}}\rangle$.
Define
$\langle\langle\alpha_{11}\rangle\rangle=\langle\alpha_{11},\alpha_{21},\dots,\alpha_{r1}\rangle$.
Therefore $\langle\langle\alpha_{11}\rangle\rangle\simeq
C_{p^{a_1}}\times C_{p^{a_2}}\times\cdots\times C_{p^{a_r}}$ is a
normal subgroup of type (4). Define
$\overline{\langle\langle\alpha_{11}\rangle\rangle}=(C_p)^t\cdot
\langle\alpha_{r+11},\dots,\alpha_{s1}\rangle$.

\emph{Sub-subcase II.3.3.2.}
$\alpha_{22}=\alpha_{31}^{p^{a_3-1}},\dots,\alpha_{r-12}=\alpha_{r1}^{p^{a_r-1}},\alpha_{r2}\notin
H^p$. Then there exist $\alpha_{r2},\dots,\alpha_{r\ell}\in H$ for
some $\ell:2\leq \ell\leq p$ such that
$[\alpha_{rj},\alpha]=\alpha_{rj+1}$, where $1\leq j\leq \ell-1$ and
$\alpha_{r\ell}\ne 1$ is central. Define
$\langle\langle\alpha_{11}\rangle\rangle=\langle\alpha_{11},\alpha_{21},\dots,\alpha_{r1},\alpha_{r2},\dots,\alpha_{r\ell}\rangle$.
Therefore $\langle\langle\alpha_{11}\rangle\rangle\simeq
C_{p^{a_1}}\times C_{p^{a_1}}\times\cdots C_{p^{a_r}}\times
(C_p)^{\ell-1}$ is a normal subgroup of type (3). Define
$\overline{\langle\langle\alpha_{11}\rangle\rangle}=(C_p)^{t-\ell+1}\cdot
\langle\alpha_{r+11},\dots,\alpha_{s1}\rangle$, where
$(C_p)^{t-\ell+1}$ is the complement of $(C_p)^{\ell-1}$ in
$(C_p)^t$.

\emph{Step III.} Put $H_1=\langle\langle\alpha_{11}\rangle\rangle$
and $H_2=\overline{\langle\langle\alpha_{11}\rangle\rangle}$. Note
that $H_1\cap H_2=\{1\}$. However, $H_2$ may not be a normal
subgroup of $G$. That is why we need to show that there exist a
commutator chain $\mathcal H_1$ and a normal subgroup $\mathcal H_2$
of $G$ such that $H=\mathcal H_1\times \mathcal H_2$. In this Step,
we will describe a somewhat algorithmic approach which replaces the
generators of $H$ until the desired result is obtained.

Assume henceforth that $H_2$ is not normal in $G$. Then there exists
a generator $\beta\in H_2$ such that $\alpha^{-1}\beta\alpha=hh_1$
for some $h\in H_2,h_1\in H_1,h_1\notin H_2$. Since $h=\beta h_2$
for some $h_2\in H_2$, we get $[\beta,\alpha]=h_1h_2$.

Let us assume first that $\ord(\beta)=p$. If $h_1\in H^p$, then
$h_2\notin [H(p),\alpha]$, otherwise $[H(p),\alpha]\cap
H^p\ne\{1\}$. In other words, $h_2$ does not appear in similar
chains, so we can simply put $h_1h_2$, instead of $h_2$, as a
generator of $H_2$. In this way we obtain a group that is
$G$-isomorphic to $H_2$. Thus we get that $[\beta,\alpha]$ is in
this new copy of $H_2$. Similarly, if $h_1\in H(p)$ and $h_2\notin
[H(p),\alpha]$, we can obtain a new copy of $H_2$ such that
$[\beta,\alpha]$ is in $H_2$. If $h_2\in [H(p),\alpha]$, we may
assume that $[\beta,\alpha]\in H_1$. In this case
$\langle\langle\alpha_{11}\rangle\rangle$ must be of type (3). Let
$\langle\langle\alpha_{11}\rangle\rangle\simeq C_{p^{a_1}}\times
C_{p^{a_2}}\times\cdots\times C_{p^{a_k}}\times(C_p)^s$, be
generated by elements
$\alpha_{11},\dots,\alpha_{k1},\alpha_{k2},\dots,\alpha_{ks+1}$ with
relations given in the statement of the Lemma. Assume that
$\alpha_{k\ell}=[\beta,\alpha]$ for some $\ell:2\leq\ell\leq s+1$.
If $\ell> 2$, replace $\beta$ with
$\beta'=\beta\alpha_{k\ell-1}^{-1}$. Hence $[\beta',\alpha]=1$. If
$\ell=2$, we can put $\alpha_{k1}'=\alpha_{k1}\beta^{-1}$, instead
of $\alpha_{k1}$, as a generator of $H_1$. In this way we obtain a
group of type (4), since $[\alpha_{k1}',\alpha]=1$. Clearly,
$[\beta,\alpha]$ is not in this new commutator chain $\mathcal H_1$.
It is not hard to see that with similar replacements we can treat
the general case
$[\beta,\alpha]=\prod_i\alpha_{i1}^{p^{a_i-1}c_i}\cdot\prod_j\alpha_{kj}$.
Thus we obtain the decomposition $H=\mathcal H_1\times \mathcal
H_2$, where $\mathcal H_1$ and $\mathcal H_2$ are normal subgroups
of $G$.

Next, we are going to assume that $\ord(\beta)>p$. According to the
definition of the commutator chain of $\alpha_{11}$ we need to
consider the three cases of Step II separately.

\emph{Case III.1.} $\alpha_{12}=\alpha_{11}^{p^{a_1-1}c_1}$ for some
$c_1:1\leq c_1\leq p-1$. Here we must have
$h_1=\alpha_{11}^{p^{a_1-1}d_1}$ for some $d_1:1\leq d_1\leq p-1$.
We can replace $\beta$ with $\beta'=\beta\alpha_{11}^{-d_1/c_1}$, so
$[\beta',\alpha]=h_2$.

\emph{Case III.2.} $\alpha_{12}\notin H^p$. If $h_1=\prod_{j\geq
2}\alpha_{1j}^{d_j}$ for some $d_j:0\leq d_j\leq p-1$, we can
replace $\beta$ with $\beta'=\beta\prod_{j\geq
2}\alpha_{1j-1}^{-d_j}$. Hence $[\beta',\alpha]=h_2$. Thus we reduce
the considerations to the case $h_1=\alpha_{11}^{p^{a_1-1}d_1}$ for
some $d_1:0\leq d_1\leq p-1$. We have now three possibilities for
$h_2$.

\emph{Subcase III.2.1.} Let $h_2\notin H^p$ and $h_2\notin
[H,\alpha]$. We can put $h_1h_2$, instead of $h_2$, as a generator
of $H_2$. In this way we obtain a group that is $G$-isomorphic to
$H_2$. Thus we get that $[\beta,\alpha]$ is in this new copy of
$H_2$.

\emph{Subcase III.2.2.} Let $h_2\notin H^p$ and $h_2\in [H,\alpha]$,
i.e., there exists $\gamma\notin H^p$ such that
$[\gamma,\alpha]=h_2$. Put $\beta'=\beta\gamma^{-1}$. Then
$[\beta',\alpha]=h_1=\alpha_{11}^{p^{a_1-1}d_1}$. Hence the
commutator chain of $\alpha_{11}$ is contained in the commutator
chain $\langle\langle\beta'\rangle\rangle$ which is a normal
subgroup of $G$ of type (3).

\emph{Subcase III.2.3.} Let $h_2\in H^p$, i.e., $h_2=\prod_{i\in
B}\alpha_{i1}^{p^{a_i-1}d_i}$, where $B=\{i:\alpha_{i1}\in
H_2\},0\leq d_i\leq p-1$. We can replace $\alpha_{11}$ with
$\alpha_{11}'=\alpha_{11}^{d_1}\prod_{i\in
B}\alpha_{i1}^{p^{a_i-a_1}d_i}$. Now we have
$[\beta,\alpha]=\alpha_{11}'^{p^{a_1-1}}$, so the commutator chain
of $\alpha_{11}'$ is contained in the commutator chain
$\langle\langle\beta\rangle\rangle$ which is a normal subgroup of
$G$ of type (3).

\emph{Case III.3.} $\alpha_{12}\in H^p$. We have that either
$\langle\langle\alpha_{11}\rangle\rangle\simeq C_{p^{a_1}}\times
C_{p^{a_2}}\times\cdots\times C_{p^{a_r}}$ is a normal subgroup of
type (4), or $\langle\langle\alpha_{11}\rangle\rangle\simeq
C_{p^{a_1}}\times C_{p^{a_1}}\times\cdots C_{p^{a_r}}\times
(C_p)^{\ell-1}$ is a normal subgroup of type (3).

Similarly to Case III.2, if $h_1$ is a product of elements of order
$p$ that are not in $\langle\alpha_{11}^{p^{a_1-1}}\rangle$, by a
suitable change of the generator $\beta$ we will obtain
$[\beta,\alpha]=h_2$. Thus we again reduce the considerations to the
case $h_1=\alpha_{11}^{p^{a_1-1}d_1}$ for some $d_1:0\leq d_1\leq
p-1$. We have three possibilities for $h_2$, which are identical to
the three subcases in Case III.2. The only slight difference is that
the new commutator chain here can be either of type (3) or (4).

In this way, we have investigated all possibilities for the proper
construction of the normal factors of $H$. The construction is
algorithmic in nature. When we define a new commutator chain
$\langle\langle\beta'\rangle\rangle$ or
$\langle\langle\beta\rangle\rangle$ (as in Subcases III.2.2 and
III.2.3), we have to start the same process all over again until we
can not get a new commutator chain that contains the previous one.
Denote by $\mathcal H_1$ the last commutator chain obtained by the
described algorithm from $H_1$. We have that $\mathcal H_1$ is a
normal subgroup of $G$ of type (1)--(4). Denote by $\mathcal H_2$
the subgroup obtained from $H_2$ by the replacements described
above. Then $H$ is a direct product of $\mathcal H_1$ and $\mathcal
H_2$, where $\mathcal H_2$ is normal in $G$. Proceeding by induction
we will obtain the decomposition given in the statement. We are
done.

\section{Proof of Theorem \ref{t1.4}}
\label{4}

If char $K=p>0$, we can apply Theorem \ref{t2.3}. Therefore, we will
assume that char $K\ne p$.

According to Lemma \ref{lemma}, $H\simeq\mathcal
H_1\times\cdots\times\mathcal H_t$, where $\mathcal
H_1,\dots,\mathcal H_t$ are normal subgroups of $G$ that are
isomorphic to any of the four types, described in Lemma \ref{lemma}.

Let $V$ be a $K$-vector space whose dual space $V^*$ is defined as
$V^*=\bigoplus_{g\in G}K\cdot x(g)$ where $G$ acts on $V^*$ by
$h\cdot x(g)=x(hg)$ for any $h,g\in G$. Thus $K(V)^G=K(x(g):g\in
G)^G=K(G)$.

Now, for any subgroup $\mathcal H_i (1\leq i\leq t)$ we can define a
faithful representation subspace $V_i=\bigoplus_{1\leq j\leq
k_i}K\cdot Y_j$, where $k_i$ is the number of the generators of
$\mathcal H_i$ as an abelian group. (For the details see Cases
I-IV.) Therefore, $\bigoplus_{1\leq i\leq t}V_i$ is a faithful
representation space of the subgroup $H$.

Next, for any  subgroup $\mathcal H_i (1\leq i\leq t)$ we define
$x_{jk}=\alpha^k\cdot Y_j$ for $1\leq j\leq k_i,0\leq k\leq p-1$.
Define $W_i=\bigoplus_{j,k}K\cdot x_{jk}\subset V^*$. Then
$W=\bigoplus_{1\leq i\leq t}W_i$ is a faithful $G$-subspace of
$V^*$. Thus, by Theorem \ref{t2.1} it suffices to show that $W^G$ is
rational over $K$. Note that
$W^G=(W^H)^{\langle\alpha\rangle}=((\dots(W^{\mathcal
H_1})^{\mathcal H_2}\dots)^{\mathcal
H_t})^{\langle\alpha\rangle}=((\dots(W_1^{\mathcal
H_1}\bigoplus_{2\leq j\leq t}W_j)^{\mathcal H_2}\dots)^{\mathcal
H_t})^{\langle\alpha\rangle}=\cdots=\bigoplus_{1\leq j\leq
t}(W_j^{\mathcal H_j})^{\langle\alpha\rangle}.$ Therefore, we need
to calculate $W_j^{\mathcal H_j}$ when $\mathcal H_j$ is isomorphic
to any of the four types, described in Lemma \ref{lemma}. Finally,
we will show that the action of $\alpha$ on $W^H$ can be linearized.

\emph{Case I.} Assume that $\mathcal H_1$ is of the type (3), i.e.,
$\mathcal H_1\simeq C_{p^{a_1}}\times C_{p^{a_2}}\times\cdots\times
C_{p^{a_k}}\times(C_p)^s$ for some $k\geq 1,a_i\geq 2,s\geq 1$.
Denote by $\alpha_1,\dots,\alpha_k$ the generators of
$C_{p^{a_1}}\times\cdots\times C_{p^{a_k}}$, and by
$\alpha_{k+1},\dots,\alpha_{k+s}$ the generators of $(C_p)^s$.
According to Lemma \ref{lemma}, we have the relations
$[\alpha_i,\alpha]=\alpha_{i+1}^{p^{a_{i+1}-1}}\in Z(G)$ for $1\leq
i\leq k-1;  [\alpha_{k+j},\alpha]=\alpha_{k+j+1}$ for $0\leq j\leq
s-1$; and $\alpha_{k+s}\in Z(G)$. Because of the frequent use of
$k+s$ in this case, we put $r=k+s$.

We divide the proof into several steps.

\emph{Step 1.} Define $X_1,X_2,\dots,X_{r}\in V^*$ by
\begin{equation*}
X_j=\sum_{\ell_1,\dots,\ell_{r}}x\left(\prod_{i\ne
j}\alpha_i^{\ell_i}\right),\quad \text{for}\ 1\leq j\leq r.
\end{equation*}

Note that $\alpha_i\cdot X_j=X_j$ for $j\ne i$. Let
$\zeta_{p^{a_i}}\in K$ be a primitive $p^{a_i}$-th root of unity for
$1\leq i\leq k$, and let $\zeta$ be a primitive $p$-th root of
unity. Define $Y_1,Y_2,\dots,Y_{r}\in V^*$ by

\begin{equation*}
Y_i=\sum_{m=0}^{p^{a_i}-1}\zeta_{p^{a_i}}^{-m}\alpha_i^m\cdot X_i,~
Y_j=\sum_{m=0}^{p-1}\zeta^{-m}\alpha_j^m\cdot X_j
\end{equation*}
for $1\leq i\leq k$ and $k+1\leq j\leq r$.

It follows that {\allowdisplaybreaks\begin{align*} \alpha_i\ :\
&Y_i\mapsto\zeta_{p^{a_i}} Y_i,~ Y_j\mapsto Y_j,\ \text{for}\ j\ne
i\
\text{and}\ 1\leq i\leq k\\
 \alpha_j\ :\
&Y_j\mapsto\zeta Y_j,~ Y_i\mapsto Y_i,\ \text{for}\ i\ne j\
\text{and}\ k+1\leq j\leq r.
\end{align*}}
Thus $V_1=\bigoplus_{1\leq j\leq r}K\cdot Y_j$ is a faithful
representation space of the subgroup $\mathcal H_1$.

Define $x_{ji}=\alpha^i\cdot Y_j$ for $1\leq j\leq r,0\leq i\leq
p-1$. Recall that $[\alpha_i,\alpha]=\alpha_{i+1}^{p^{a_{i+1}-1}}\in
Z(G)$ for $1\leq i\leq k-1;  [\alpha_{k+j},\alpha]=\alpha_{k+j+1}$
for $0\leq j\leq s-1$; and $\alpha_{r}\in Z(G)$. Hence
$$\alpha^{-i}\alpha_j\alpha^i=\alpha_j\alpha_{j+1}^{ip^{a_{i+1}-1}},\quad \text{for}\ 1\leq j\leq k-1,1\leq i\leq p-1$$
and
$$\alpha^{-i}\alpha_j\alpha^i=\alpha_j\alpha_{j+1}^{\binom{i}{1}}\alpha_{j+2}^{\binom{i}{2}}\cdots
\alpha_{r}^{\binom{i}{r-j}},\quad \text{for}\ k\leq j\leq r-1,1\leq
i\leq p-1.$$

It follows that {\allowdisplaybreaks\begin{align*}\alpha_{\ell}\ :\
&x_{\ell i}\mapsto\zeta_{p^{a_\ell}} x_{\ell i},~x_{\ell+1
i}\mapsto\zeta^i x_{\ell+1i},~ x_{ji}\mapsto x_{ji},\ \text{for}\
1\leq\ell\leq k-1\
\text{and}\ j\ne\ell,\ell+1,\\
\alpha_k\ :\ &x_{ki}\mapsto \zeta_{p^{a_k}} x_{ki},~ x_{wi}\mapsto
\zeta^{\binom{i}{w-k}} x_{wi},~ x_{v
i}\mapsto x_{v i},\ \text{for}\ 1\leq v\leq k-1,k+1\leq w\leq r,\\
\alpha_m\ :\ &x_{ui}\mapsto \zeta^{\binom{i}{u-m}} x_{ui},~ x_{v
i}\mapsto x_{v i},\ \text{for}\ k+1\leq m\leq r,1\leq v\leq m-1,m\leq u\leq r,\\
\alpha\ :\ &x_{j0}\mapsto x_{j1}\mapsto\cdots\mapsto x_{jp-1}\mapsto
\zeta_{p^{c_j}}^{b_j}x_{j0},\ \text{for}\ 1\leq j\leq r,
\end{align*}}
where $0\leq i\leq p-1$, and $c_j,b_j$ are some integers such that
$0\leq b_j< p^{c_j}\leq p^{a_j}$.

Define $W_1=\bigoplus_{j,i}K\cdot x_{ji}\subset V^*$. As we noted in
the beginning of the proof, we need to find $W_1^{\mathcal H_1}$.

\emph{Step 2.} For $1\leq j\leq r$ and for $1\leq i\leq p-1$ define
$y_{ji}=x_{ji}/x_{ji-1}$. Thus $W_1=K(x_{j0},y_{ji}:1\leq j\leq
r,1\leq i\leq p-1)$ and for every $g\in G$
\begin{equation*}
g\cdot x_{j0}\in K(y_{ji}:1\leq j\leq r,1\leq i\leq p-1)\cdot
x_{j0},\ \text{for}\ 1\leq j\leq r
\end{equation*}
while the subfield $K(y_{ji}:1\leq j\leq r,1\leq i\leq p-1)$ is
invariant by the action of $G$, i.e.,
{\allowdisplaybreaks\begin{align*}\alpha_{\ell}\ :\ &y_{\ell+1
i}\mapsto\zeta y_{\ell+1i},~ y_{ji}\mapsto y_{ji},\ \text{for}\
1\leq\ell\leq k-1\
\text{and}\ j\ne\ell+1,\\
\alpha_m\ :\ &y_{ui}\mapsto \zeta^{\binom{i-1}{u-m-1}} y_{ui},~ y_{v
i}\mapsto y_{v i},\ \text{for}\ k\leq m\leq r-1,1\leq v\leq m,m+1\leq u\leq r,\\
\alpha_{r}\ :\ &y_{v i}\mapsto y_{v i},\ \text{for}\ 1\leq v\leq r,\\
\alpha\ :\ &y_{j1}\mapsto y_{j2}\mapsto\cdots\mapsto y_{jp-1}\mapsto
\zeta_{p^{c_j}}^{b_j}(y_{j1}\cdots y_{jp-1})^{-1},\ \text{for}\
1\leq j\leq r,
\end{align*}}
From Theorem \ref{t2.2} it follows that if $K(y_{ji}:1\leq j\leq
r,1\leq i\leq p-1)^{G}$ is rational over $K$, so is
$K(x_{j0},y_{ji}:1\leq j\leq r,1\leq i\leq p-1)^{G}$ over $K$.

Since $K$ contains a primitive $p^e$-th root of unity $\zeta_{p^e}$
where $p^e$ is the exponent of $G$, $K$ contains as well a primitive
$p^{{c_j}+1}$-th root of unity, and we may replace the variables
$y_{ji}$ by $y_{ji}/\zeta_{p^{{c_j}+1}}^{b_j}$ so that we obtain a
more convenient action of $\alpha$ without changing the actions of
$\alpha_j$'s. Namely we may assume that
\begin{align*}
\alpha\ :\ &y_{j1}\mapsto y_{j2}\mapsto\cdots\mapsto y_{jp-1}\mapsto
(y_{j1}y_{j2}\dots y_{jp-1})^{-1}\ \text{for}\ 1\leq j\leq r.
\end{align*}

Define $u_{r1}=y_{r1}^p,u_{ri}=y_{ri}/y_{ri-1}$ for $2\leq i\leq
p-1$. Then $K(y_{ji},u_{ri}:1\leq j\leq r-1,1\leq i\leq
p-1)=K(y_{ji}:1\leq j\leq r,1\leq i\leq
p-1)^{\langle\alpha_{r-1}\rangle}$. From Theorem \ref{t2.2} it
follows that if $K(y_{ji},u_{ri}:1\leq j\leq r-1,2\leq i\leq
p-1)^{G}$ is rational over $K$, so is $K(y_{ji},u_{ri}:1\leq j\leq
r-1,1\leq i\leq p-1)^{G}$ over $K$. We have the following actions
{\allowdisplaybreaks\begin{align*}\alpha_{\ell}\ :\ &u_{ri}\mapsto
u_{ri},\ \text{for}\
1\leq\ell\leq k-1,\\
\alpha_m\ :\ &u_{ri}\mapsto \zeta^{\binom{i-2}{r-m-2}} u_{ri},\
\text{for}\ 2\leq i\leq p-1\ \text{and}\ k\leq m\leq r-2,\\
\alpha\ :\ &u_{r2}\mapsto u_{r3}\mapsto\cdots\mapsto u_{rp-1}\mapsto
(u_{r1}u_{r2}^{p-1}u_{r3}^{p-2}\cdots u_{rp-1}^2)^{-1}\mapsto
u_{r1}u_{r2}^{p-2}u_{r3}^{p-3}\cdots u_{rp-2}^2u_{rp-1}.
\end{align*}}
For $2\leq i\leq p-1$ define
$$v_{ri}=u_{ri}y_{r-1i}^{-1}y_{r-2i}y_{r-3i}^{-1}\cdots
y_{k+2i}^{(-1)^{r-k}}y_{k+1i}^{(-1)^{r-k+1}},$$ and put
$v_{r1}=u_{r1}$.

With the aid of the well known property
$\binom{n}{m}-\binom{n-1}{m}=\binom{n-1}{m-1}$, it is not hard to
verify the identity {\allowdisplaybreaks\begin{align*} \
&\binom{i-2}{r-m-2}-\binom{i-1}{r-m-2}+\binom{i-1}{r-m-3}-\binom{i-1}{r-m-4}+\cdots\\
&\cdots+(-1)^{r-m-1}\binom{i-1}{2}+(-1)^{r-m}\binom{i-1}{1}+(-1)^{r-m+1}\binom{i-1}{0}=0.
\end{align*}}
It follows that {\allowdisplaybreaks\begin{align*} \alpha_m\ :\
&v_{ri}\mapsto v_{ri},\
\text{for}\ 1\leq i\leq p-1\ \text{and}\ 1\leq m\leq r-2,\\
\alpha\ :\ &v_{r2}\mapsto v_{r3}\mapsto\cdots\mapsto v_{rp-1}\mapsto
A_r\cdot(v_{r1}v_{r2}^{p-1}v_{r3}^{p-2}\cdots v_{rp-1}^2)^{-1}.
\end{align*}}
where $A_r$ is some monomial in $y_{ji}$ for $2\leq j\leq r-1,1\leq
i\leq p-1$.

Define $u_{r-11}=y_{r-11}^p,u_{r-1i}=y_{r-1i}/y_{r-1i-1}$ for $2\leq
i\leq p-1$. Then $K(y_{ji},u_{r-1i}:1\leq j\leq r-2,1\leq i\leq
p-1)=K(y_{ji}:1\leq j\leq r-1,1\leq i\leq
p-1)^{\langle\alpha_{r-2}\rangle}$. From Theorem \ref{t2.2} it
follows that if $K(y_{ji},u_{r-1i}:1\leq j\leq r-2,2\leq i\leq
p-1)^{G}$ is rational over $K$, so is $K(y_{ji},u_{r-1i}:1\leq j\leq
r-2,1\leq i\leq p-1)^{G}$ over $K$. Similarly to the definition of
$v_{ri}$, we can define $v_{r-1i}$ so that
$\alpha_m(v_{r-1i})=v_{r-1i}$ for $2\leq i\leq p-1$ and $1\leq m\leq
r-3$. It is obvious that we can proceed in the same way defining
elements $v_{r-2i}, v_{r-3i},\dots,v_{k+1i}$ such that $\alpha_m$
acts trivially on all $v_{ji}$'s for $k\leq m\leq r-3$.

Recall that the actions of $\alpha_{\ell}$ on $y_{ji}$'s for
$1\leq\ell\leq k-1$ are
\begin{align*}\alpha_{\ell}\ :\ &y_{\ell+1
i}\mapsto\zeta y_{\ell+1i},~ y_{ji}\mapsto y_{ji},\ \text{for}\
1\leq i\leq p-1, 1\leq\ell\leq k-1\ \text{and}\ j\ne\ell+1.
\end{align*}

For any $1\leq\ell\leq k-1$ define
$v_{\ell+11}=y_{\ell+11}^p,v_{\ell+1i}=y_{\ell+1i}/y_{\ell+1i-1}$,
where $2\leq i\leq p-1$. Put also $v_{1i}=y_{1i}$ for $1\leq i\leq
p-1$. Then $K(v_{ji}:1\leq j\leq r,1\leq i\leq p-1)=K(y_{ji}:1\leq
j\leq r,1\leq i\leq p-1)^{\mathcal H_1}$.

The action of $\alpha$ is given by
{\allowdisplaybreaks\begin{align*} \alpha\ :\ &v_{11}\mapsto
v_{12}\mapsto\cdots\mapsto v_{1p-1}\mapsto (v_{11}v_{12}\cdots
v_{1p-1})^{-1},\\
\ &v_{m1}\mapsto v_{m1}v_{m2}^p,~ v_{m2}\mapsto
v_{m3}\mapsto\cdots\mapsto v_{mp-1}\mapsto
A_m\cdot(v_{m1}v_{m2}^{p-1}v_{m3}^{p-2}\cdots v_{mp-1}^2)^{-1},\\
\ &\text{for}\ 2\leq m\leq r,
\end{align*}}
where $A_m$ is some monomial in $v_{k+1i},\dots,v_{m-1i}$ for
$k+2\leq m\leq r$ and $A_2=A_3=\cdots=A_{k+1}=1$. From Lemmas
\ref{l2.7} and \ref{l2.8} it follows that the action of $\alpha$ on
$K(v_{ji}:1\leq j\leq r,1\leq i\leq p-1)$ can be linearized.

\emph{Case II.} Assume that $\mathcal H_1$ is of the type (1), i.e.,
$\mathcal H_1\simeq (C_p)^{s+1}$ for some $s\geq 0$. Denote by
$\beta_1,\dots,\beta_{s+1}$ the generators of $(C_p)^{s+1}$.
According to Lemma \ref{lemma}, we have the relations
$[\beta_j,\alpha]=\beta_{j+1}$ for $1\leq j\leq s$; and
$\beta_{s+1}\in Z(G)$.

Define $X_1,X_2,\dots,X_{s+1}\in V^*$ by
\begin{equation*}
X_j=\sum_{\ell_1,\dots,\ell_{s+1}}x\left(\prod_{m\ne
j}\beta_m^{\ell_m}\right),
\end{equation*}
for $1\leq j\leq s+1$. Note that $\beta_j\cdot X_i=X_i$ for $j\ne
i$. Let $\zeta$ be a primitive $p$-th root of unity. Define
$Y_1,Y_2,\dots,Y_{s+1}\in V^*$ by
\begin{equation*}
Y_j=\sum_{r=0}^{p-1}\zeta^{-r}\beta_j^r\cdot X_j
\end{equation*}
for $1\leq j\leq s+1$.

It follows that {\allowdisplaybreaks\begin{align*} &\beta_j\ :\
Y_j\mapsto\zeta Y_j,~ Y_i\mapsto Y_i,\ \text{for}\ i\ne j\
\text{and}\ 1\leq j\leq s+1.
\end{align*}}
Thus $V_1=\bigoplus_{1\leq j\leq s+1}K\cdot Y_j$ is a representation
space of the subgroup $\mathcal H_1$.

Define $x_{ji}=\alpha^i\cdot Y_j$ for $1\leq j\leq s+1,0\leq i\leq
p-1$. Recall that $[\beta_j,\alpha]=\beta_{j-1}$. Hence
$$\alpha^{-i}\beta_j\alpha^i=\beta_j\beta_{j+1}^{\binom{i}{1}}\beta_{j+2}^{\binom{i}{2}}\cdots
\beta_{s+1}^{\binom{i}{s+1-j}}.$$

It follows that {\allowdisplaybreaks\begin{align*}\beta_1\ :\
&x_{1i}\mapsto\zeta x_{1i},~ x_{ji}\mapsto \zeta^{\binom{i}{j-1}}
x_{ji},\ \text{for}\ 2\leq j\leq s+1\
\text{and}\ 0\leq i\leq p-1,\\
\beta_j\ :\ &x_{\ell i}\mapsto x_{\ell i},~ x_{mi}\mapsto
\zeta^{\binom{i}{m-j}} x_{mi},\ \text{for}\ 1\leq \ell\leq j-1,j\leq
m\leq s+1\ \text{and}\ 0\leq i\leq p-1,\\
\alpha\ :\ &x_{j0}\mapsto x_{j1}\mapsto\cdots\mapsto x_{jp-1}\mapsto
\zeta^{b_j}x_{j0},\ \text{for}\ 1\leq j\leq s+1, 0\leq b_j\leq p-1.
\end{align*}}
Compare the actions of $\alpha,\beta_1,\dots,\beta_{s+1}$ with the
actions of $\alpha,\alpha_k,\dots,\alpha_{k+s}$ from Case I, Step 1.
They are almost the same. Apply the proof of Case I.

\emph{Case III.} Assume that $\mathcal H_1$ is of the type (2),
i.e., $\mathcal H_1\simeq C_{p^a}$ for some $a\geq 1$. Denote by
$\beta$ the generator of $C_{p^a}$. Then
$[\beta,\alpha]=\beta^{bp^{a-1}}$ for some $b:0\leq b\leq p-1$. Let
$\zeta_{p^a}\in K$ be a primitive $p^a$-th root of unity, and let
$\zeta$ be a primitive $p$-th root of unity. Define
$X=\sum_i\zeta_{p^a}^{-i}x(\beta^i)$. Then $\beta(X)=\zeta_{p^a}X$,
and define $x_i=\alpha^i\cdot X$ for $0\leq i\leq p-1$. It follows
that {\allowdisplaybreaks\begin{align*}\beta\ :\ &x_i\mapsto
\zeta_{p^a}\zeta^{ib} x_i,\ \text{for}\ 0\leq i\leq p-1,\\
\alpha\ :\ &x_0\mapsto x_1\mapsto\cdots\mapsto x_{p-1}\mapsto
\zeta_{p^a}^{c}x_0,\ \text{for}\  0\leq c\leq p^a-1.
\end{align*}}
Define $W_1=\bigoplus_iK\cdot x_i\subset V^*$. For $1\leq i\leq p-1$
define $y_i=x_i/x_{i-1}$. Thus $W_1=K(x_0,y_i:1\leq i\leq p-1)$ and
for every $g\in G$
\begin{equation*}
g\cdot x_0\in K(y_i:1\leq i\leq p-1)\cdot x_0,
\end{equation*}
while the subfield $K(y_i:1\leq i\leq p-1)$ is invariant by the
action of $G$, i.e., {\allowdisplaybreaks\begin{align*}\beta\ :\
&y_i\mapsto
\zeta^b y_i,\ \text{for}\ 1\leq i\leq p-1,\\
\alpha\ :\ &y_1\mapsto y_2\mapsto\cdots\mapsto
\zeta_{p^a}^{c}(y_1\cdots y_{p-1})^{-1},\ \text{for}\  0\leq c\leq
p^a-1.
\end{align*}}
From Theorem \ref{t2.2} it follows that if $K(y_i:1\leq i\leq
p-1)^G$ is rational over $K$, so is $K(x_0,y_i:1\leq i\leq p-1)^G$
over $K$.

Since $K$ contains a primitive $p^e$-th root of unity $\zeta_{p^e}$
where $p^e$ is the exponent of $G$, $K$ contains as well
$\zeta_{p^{a+1}}^{c}$. We may replace the variables $y_i$ by
$y_i/\zeta_{p^{a+1}}^c$ so that we obtain
\begin{align*}
\alpha\ :\ &y_1\mapsto y_2\mapsto\cdots\mapsto y_{p-1}\mapsto
(y_{1}y_{2}\dots y_{p-1})^{-1}.
\end{align*}

Define $u_1=y_1^p,u_i=y_i/y_{i-1}$ for $2\leq i\leq p-1$. Then
$K(u_i:1\leq i\leq p-1)=K(y_i:1\leq i\leq
p-1)^{\langle\beta\rangle}$. The action of $\alpha$ is given by
{\allowdisplaybreaks\begin{align*} \alpha\ :\ &u_1\mapsto u_1u_2^p,~
u_{2}\mapsto u_{3}\mapsto\cdots\mapsto u_{p-1}\mapsto
(u_{1}u_{2}^{p-1}u_{3}^{p-2}\cdots u_{p-1}^2)^{-1},
\end{align*}}
From Lemma \ref{l2.7} (or \ref{l2.8}) it follows that the action of
$\alpha$ can be linearized.

\emph{Case IV.} Assume that $\mathcal H_1$ is of the type (4), i.e.,
$\mathcal H_1\simeq C_{p^{a_1}}\times C_{p^{a_2}}\times\cdots\times
C_{p^{a_k}}$ for some $k\geq 2$. Denote by $\alpha_1,\dots,\alpha_k$
the generators of $\mathcal H_1$. According to Lemma \ref{lemma}, we
have the relations
$[\alpha_i,\alpha]=\alpha_{i+1}^{p^{a_{i+1}-1}}\in Z(G)$ for $1\leq
i\leq k-1;
[\alpha_k,\alpha]=\prod_{j=1}^{k}\alpha_j^{p^{a_j-1}c_j}\in Z(G)$
for some $0\leq c_j\leq p-1$.

Similarly to the previous cases, define $Y_1,Y_2,\dots,Y_{k}\in V^*$
so that {\allowdisplaybreaks\begin{align*} \alpha_i\ :\
&Y_i\mapsto\zeta_{p^{a_i}} Y_i,~ Y_j\mapsto Y_j,\ \text{for}\ j\ne
i\ \text{and}\ 1\leq i\leq k.
\end{align*}}
Thus $V_1=\bigoplus_{1\leq j\leq k}K\cdot Y_j$ is a faithful
representation space of the subgroup $\mathcal H_1$.

Next, define $x_{ji}=\alpha^i\cdot Y_j$ for $1\leq j\leq k,0\leq
i\leq p-1$. Note that
$$\alpha^{-i}\alpha_j\alpha^i=\alpha_j\alpha_{j+1}^{ip^{a_{j+1}-1}},\quad \text{for}\ 1\leq j\leq k-1,1\leq i\leq p-1$$
and
$$\alpha^{-i}\alpha_k\alpha^i=\alpha_k\prod_{j=1}^{k}\alpha_j^{ip^{a_j-1}c_j},\quad \text{for}\ 1\leq
i\leq p-1.$$

It follows that {\allowdisplaybreaks\begin{align*}\alpha_{\ell}\ :\
&x_{\ell i}\mapsto\zeta_{p^{a_\ell}} x_{\ell i},~x_{\ell+1
i}\mapsto\zeta^i x_{\ell+1i},~ x_{ji}\mapsto x_{ji},\ \text{for}\
1\leq\ell\leq k-1\
\text{and}\ j\ne\ell,\ell+1,\\
\alpha_k\ :\ &x_{ki}\mapsto \zeta_{p^{a_k}}\zeta^{ic_k} x_{ki},~
x_{ji}\mapsto \zeta^{ic_j} x_{ji},\ \text{for}\ 1\leq j\leq k-1,\\
\alpha\ :\ &x_{j0}\mapsto x_{j1}\mapsto\cdots\mapsto x_{jp-1}\mapsto
\zeta_{p^{a_j}}^{b_j}x_{j0},\ \text{for}\ 1\leq j\leq k,
\end{align*}}
where $0\leq i\leq p-1,0\leq c_j\leq p-1$ and $0\leq b_j\leq
p^{a_j}-1$.

Define $W_1=\bigoplus_{j,i}K\cdot x_{ji}\subset V^*$, and for $1\leq
i\leq p-1$ define $y_i=x_i/x_{i-1}$. Thus $W_1=K(x_{j0},y_{ji}:1\leq
j\leq k,1\leq i\leq p-1)$ and for every $g\in G$
\begin{equation*}
g\cdot x_{j0}\in K(y_{ji}:1\leq j\leq k,1\leq i\leq p-1)\cdot
x_{j0},\ \text{for}\ 1\leq j\leq k
\end{equation*}
while the subfield $K(y_{ji}:1\leq j\leq k,1\leq i\leq p-1)$ is
invariant by the action of $G$, i.e.,
{\allowdisplaybreaks\begin{align*}\alpha_{\ell}\ :\ &y_{\ell+1
i}\mapsto\zeta y_{\ell+1i},~ y_{ji}\mapsto y_{ji},\ \text{for}\
1\leq i\leq p-1,~ 1\leq\ell\leq k-1\
\text{and}\ j\ne\ell+1,\\
\alpha_k\ :\ &y_{ji}\mapsto \zeta^{c_j} y_{ji},\ \text{for}\ 1\leq
i\leq p-1,~ 1\leq j\leq k,\\
\alpha\ :\ &y_{j1}\mapsto y_{j2}\mapsto\cdots\mapsto y_{jp-1}\mapsto
\zeta_{p^{a_j}}^{b_j}(y_{j1}\cdots y_{jp-1})^{-1}
\end{align*}}
From Theorem \ref{t2.2} it follows that if $K(y_{ji}:1\leq j\leq
k,1\leq i\leq p-1)^{G}$ is rational over $K$, so is
$K(x_{j0},y_{ji}:1\leq j\leq k,1\leq i\leq p-1)^{G}$ over $K$. As
before, we can again assume that $\alpha$ acts in this way:
\begin{align*}
\alpha\ :\ &y_{j1}\mapsto y_{j2}\mapsto\cdots\mapsto y_{jp-1}\mapsto
(y_{j1}y_{j2}\dots y_{jp-1})^{-1}.
\end{align*}

Now, assume that $0<c_1\leq p-1$. For $2\leq j\leq k$ choose $e_j$
such that $c_1e_j+c_j\equiv 0\pmod p$, and define $u_{1i}=y_{1i},
u_{ji}=y_{1i}^{e_j}y_{ji}$. It follows that
{\allowdisplaybreaks\begin{align*}\alpha_{\ell}\ :\ &u_{\ell+1
i}\mapsto\zeta u_{\ell+1i},~ u_{ji}\mapsto u_{ji},\ \text{for}\
1\leq i\leq p-1,~ 1\leq\ell\leq k-1\
\text{and}\ j\ne\ell+1,\\
\alpha_k\ :\ &u_{1i}\mapsto\zeta^{c_1} u_{1i}, u_{ji}\mapsto
u_{ji},\ \text{for}\ 1\leq i\leq p-1,~ 2\leq j\leq k.
\end{align*}}
Define $v_{j1}=u_{j1}^p,v_{ji}=u_{ji}/u_{ji-1}$ for $2\leq i\leq
p-1,1\leq j\leq k$. Then $K(v_{ji}:1\leq j\leq k,1\leq i\leq
p-1)=K(u_{ji}:1\leq j\leq k,1\leq i\leq p-1)^{\mathcal H_1}$. The
action of $\alpha$ is given by {\allowdisplaybreaks\begin{align*}
\alpha\ :\ &v_{j1}\mapsto v_{j1}v_{j2}^p,~ v_{j2}\mapsto
v_{j3}\mapsto\cdots\mapsto v_{jp-1}\mapsto
(v_{j1}v_{j2}^{p-1}v_{j3}^{p-2}\cdots v_{jp-1}^2)^{-1},\\
\ &\text{for}\ 2\leq j\leq k.
\end{align*}}
From Lemma \ref{l2.8} it follows that the action of $\alpha$ on
$K(v_{ji}:1\leq j\leq k,1\leq i\leq p-1)$ can be linearized.

Finally, let $c_1=0$. Define
$v_{j1}=u_{j1}^p,v_{ji}=u_{ji}/u_{ji-1}$ for $2\leq i\leq p-1,2\leq
j\leq k$. Then $K(u_{1i},v_{ji}:2\leq j\leq k,1\leq i\leq
p-1)=K(u_{ji}:1\leq j\leq k,1\leq i\leq p-1)^{\mathcal H_1}$. The
action of $\alpha$ again can be linearized as before. We are done.

\section{Proof of Theorem \ref{t1.5}}
\label{5}

By studying the classification of all groups of order $p^5$ made by
James in \cite{Ja}, we see that the non abelian groups with an
abelian subgroup of index $p$ and that are not direct products of
smaller groups are precisely the groups from the isoclinic families
with numbers $2,3,4$ and $9$. Notice that all these groups satisfy
the conditions of Theorem \ref{t1.4}. The isoclinic family $8$
contains only the group $\Phi_8(32)$ which is metacyclic, so we can
apply Theorem \ref{t1.2}. It is not hard to see that there are no
other groups of order $p^5$, containing a normal abelian subgroup
$H$ such that $G/H$ is cyclic.

The groups of order $p^6$ with an abelian subgroup of index $p$ and
that are not direct products of smaller groups are precisely the
groups from the isoclinic families with numbers $2,3,4$ and $9$.
Again, all these groups satisfy the conditions of Theorem
\ref{t1.4}. The groups of order $p^6$, containing a normal abelian
subgroup $H$ such that $G/H$ is cyclic of order $>p$ are precisely
the groups from the isoclinic families with numbers $8$ and $14$.
Note that the groups $\Phi_8(42),\Phi_8(33),\Phi_{14}(42)$ are
metacyclic, and the group $\Phi_8(321)a$ is a direct product of the
metacyclic group $\Phi_8(32)$ and the cyclic group $C_p$. Therefore,
we need to consider the remaining groups, whose presentations we
write down for convenience of the reader. {\allowdisplaybreaks
\begin{align*}
\Phi_8(321)b=& \langle \alpha_1,\alpha_2,\beta,\gamma:[\alpha_1,\alpha_2]=\beta=\alpha_1^p,~[\beta,\alpha_2]=\beta^p=\gamma^p,~\alpha_2^{p^2}=\beta^{p^2}=1\rangle, \\
\Phi_8(321)c_r=& \langle \alpha_1,\alpha_2,\beta:[\alpha_1,\alpha_2]=\beta,~[\beta,\alpha_2]^{r+1}=\beta^{p(r+1)}=\alpha_1^{p^2},~\alpha_2^{p^2}=\beta^{p^2}=1\rangle, \\
\Phi_8(321)c_{p-1}=& \langle \alpha_1,\alpha_2,\beta:[\alpha_1,\alpha_2]=\beta,~[\beta,\alpha_2]=\beta^p=\alpha_2^{p^2},~\alpha_1^{p^2}=\beta^{p^2}=1\rangle, \\
\Phi_8(222)=& \langle
\alpha_1,\alpha_2,\beta:[\alpha_1,\alpha_2]=\beta,~[\beta,\alpha_2]=\beta^p,~\alpha_1^{p^2}=\alpha_2^{p^2}=\beta^{p^2}=1\rangle,\\
\Phi_{14}(321)=& \langle \alpha_1,\alpha_2,\beta:[\alpha_1,\alpha_2]=\beta,~\alpha_1^{p^2}=\beta^p,~\alpha_2^{p^2}=\beta^{p^2}=1\rangle, \\
\Phi_{14}(222)=& \langle
\alpha_1,\alpha_2,\beta:[\alpha_1,\alpha_2]=\beta,~\alpha_1^{p^2}=\alpha_2^{p^2}=\beta^{p^2}=1\rangle.
\end{align*}}

\emph{Case I.} $G=\Phi_8(321)b$. Denote by $H$ the abelian normal
subgroup of $G$ generated by $\alpha_1$ and $\gamma$. Then
$H=\langle\alpha_1,\gamma\beta^{-1}\rangle\simeq C_{p^3}\times C_p$
and $G/H=\langle\alpha_2\rangle\simeq C_{p^2}$.

Let $V$ be a $K$-vector space whose dual space $V^*$ is defined as
$V^*=\bigoplus_{g\in G}K\cdot x(g)$ where $G$ acts on $V^*$ by
$h\cdot x(g)=x(hg)$ for any $h,g\in G$. Thus $K(V)^{G}=K(x(g):g\in
G)^{G}=K(G)$.

Define $X_1,X_2\in V^*$ by
\begin{equation*}
X_1=\sum_{i=0}^{p-1}x((\gamma\beta^{-1})^i),~
X_2=\sum_{i=0}^{p^3-1}x(\alpha_1^i).
\end{equation*}
Note that $\gamma\beta^{-1}\cdot X_1=X_1$ and $\alpha_1\cdot
X_2=X_2$.

Let $\zeta_{p^3}\in K$ be a primitive $p^3$-th root of unity and put
$\zeta=\zeta_{p^3}^{p^2}$, a primitive $p$-th root of unity. Define
$Y_1,Y_2\in V^*$ by
\begin{equation*}
Y_1=\sum_{i=0}^{p^3-1}\zeta_{p^3}^{-i}\alpha_1^i\cdot X_1,~
Y_2=\sum_{i=0}^{p-1}\zeta^{-i}(\gamma\beta^{-1})^i\cdot X_2.
\end{equation*}

It follows that {\allowdisplaybreaks \begin{align*}
\alpha_1\ :\ &Y_1\mapsto\zeta_{p^3} Y_1,~Y_2\mapsto Y_2,\\
\gamma\beta^{-1}\ :\ &Y_1\mapsto Y_1,~Y_2\mapsto\zeta Y_2,\\
\gamma\ :\ &Y_1\mapsto\zeta_{p^2} Y_1,~Y_2\mapsto\zeta Y_2.
\end{align*}}
Thus $K\cdot Y_1+K\cdot Y_2$ is a representation space of the
subgroup $H$.

Define $x_i=\alpha_2^i\cdot Y_1,y_i=\alpha_2^i\cdot Y_2$ for $0\leq
i\leq p^2-1$. From the relations
$\alpha_1\alpha_2^i=\alpha_2^i\alpha_1\beta^i\beta^{\binom{i}{2}p}$
it follows that {\allowdisplaybreaks \begin{align*}
\alpha_1\ :\ &x_i\mapsto\zeta_{p^3}\zeta_{p^2}^i\zeta^{\binom{i}{2}} x_i,~ y_i\mapsto y_i\\
\gamma\ :\ &x_i\mapsto\zeta_{p^2} x_i,~ y_i\mapsto\zeta y_i,\\
\alpha_2\ :\ &x_0\mapsto x_1\mapsto\cdots\mapsto x_{p^2-1}\mapsto x_0,\\
&y_0\mapsto y_1\mapsto\cdots\mapsto y_{p^2-1}\mapsto y_0,
\end{align*}}
for $0\leq i\leq p^2-1$.

We find that $Y=(\bigoplus_{0\leq i\leq p^2-1}K\cdot
x_i)\oplus(\bigoplus_{0\leq i\leq p^2-1}K\cdot y_i)$ is a faithful
$G$-subspace of $V^*$. Thus, by Theorem \ref{t2.1}, it suffices to
show that $K(x_i,y_i:0\leq i\leq p^2-1)^{G}$ is rational over $K$.

For $1\leq i\leq p^2-1$, define $u_i=x_i/x_{i-1}$ and
$v_i=y_i/y_{i-1}$. Thus $K(x_i,y_i:0\leq i\leq
p^2-1)=K(x_0,y_0,u_i,v_i:1\leq i\leq p^2-1)$ and for every $g\in G$
\begin{equation*}
g\cdot x_0\in K(u_i,v_i:1\leq i\leq p^2-1)\cdot x_0,~ g\cdot y_0\in
K(u_i,v_i:1\leq i\leq p^2-1)\cdot y_0,
\end{equation*}
while the subfield $K(u_i,v_i:1\leq i\leq p^2-1)$ is invariant by
the action of $G$. Thus $K(x_i,y_i:0\leq i\leq
p^2-1)^{G}=K(u_i,v_i:1\leq i\leq p^2-1)^{G}(u,v)$ for some $u,v$
such that $\alpha_1(v)=\gamma(v)=\alpha_2(v)=v$ and
$\alpha_1(u)=\gamma(u)=\alpha_2(u)=u$. We have now
{\allowdisplaybreaks
\begin{align*}
\alpha_1\ :\ &u_i\mapsto\zeta_{p^2}\zeta^{i-1} u_i,~ v_i\mapsto v_i,\\
\tag{5.1} \gamma\ :\ &u_i\mapsto u_i,~ v_i\mapsto v_i,\\
\alpha_2\ :\ &u_1\mapsto u_2\mapsto\cdots\mapsto u_{p^2-1}\mapsto (u_1u_2\cdots u_{p^2-1})^{-1},\\
&v_1\mapsto v_2\mapsto\cdots\mapsto v_{p^2-1}\mapsto (v_1v_2\cdots
v_{p^2-1})^{-1},
\end{align*}}
for $1\leq i\leq p^2-1$. From Theorem \ref{t2.2} it follows that if
$K(u_i,v_i:1\leq i\leq p^2-1)^{G}(u,v)$ is rational over $K$, so is
$K(x_i,y_i:0\leq i\leq p^2-1)^{G}$ over $K$.

Since $\gamma$ acts trivially on $K(u_i,v_i:1\leq i\leq p^2-1)$, we
find that $K(u_i,v_i:1\leq i\leq p^2-1)^{G}=K(u_i,v_i:1\leq i\leq
p^2-1)^{\langle\alpha_1,\alpha_2\rangle}$.

Now, consider the metacyclic $p$-group $\widetilde
G=\langle\sigma,\tau:\sigma^{p^3}=\tau^{p^2}=1,\tau^{-1}\sigma\tau=\sigma^{k},k=1+p\rangle$.

Define $X=\sum_{0\leq j\leq
p^3-1}\zeta_{p^3}^{-j}x(\sigma^j),V_i=\tau^i X$ for $0\leq i\leq
p^2-1$. It follows that
\begin{eqnarray*}
\sigma&:&V_i\mapsto \zeta_{p^3}^{k^i}V_i,\\
\tau&:&V_0\mapsto V_1\mapsto\cdots\mapsto V_{p^2-1}\mapsto V_0.
\end{eqnarray*}
Note that $K(V_0,V_1,\dots,V_{p^2-1})^{\widetilde G}$ is rational by
Theorem \ref{t2.6}.

Define $U_i=V_i/V_{i-1}$ for $1\leq i\leq p^2-1$. Then
$K(V_0,V_1,\dots,V_{p^2-1})^{\widetilde G}=K(U_1,U_2,\dots,$
$U_{p^2-1})^{\widetilde G}(U)$ where
\begin{eqnarray*}
\sigma&:&U_i\mapsto \zeta_{p^3}^{k^i-k^{i-1}}U_i,~ U\mapsto U\\
\tau&:&U_1\mapsto U_2\mapsto\cdots\mapsto U_{p^2-1}\mapsto
(U_1U_2\cdots U_{p^2-1})^{-1},~ U\mapsto U.
\end{eqnarray*}

Notice that $k^i-k^{i-1}=(1+p)^{i-1}p\equiv (1+(i-1)p)p\pmod{p^3}$,
so $\zeta_{p^3}^{k^i-k^{i-1}}=\zeta_{p^2}^{1+(i-1)p}$. Compare
Formula (5.1) (i.e., the actions of $\alpha_1,\alpha_2$ on
$K(u_i:1\leq i\leq p^2-1)$) with the actions of $\widetilde G$ on
$K(U_i:1\leq i\leq p^2-1)$. They are the same. Hence, according to
Theorem \ref{t2.6}, we get that $K(u_1,\dots,u_{p^2-1})^{G}(u)\cong
K(U_1,\dots,U_{p^2-1})^{\widetilde
G}(U)=K(V_0,V_1,\dots,V_{p^2-1})^{\widetilde G}$ is rational over
$K$. Since by Lemma \ref{l2.7} we can linearize the action of
$\alpha_2$ on $K(v_i:1\leq i\leq p^2-1)$, we obtain finally that
$K(u_i,v_i:1\leq i\leq p^2-1)^{\langle\alpha_1,\alpha_2\rangle}$ is
rational over $K$.

\emph{Case II.} $G=\Phi_8(321)c_r$. Denote by $H$ the abelian normal
subgroup of $G$ generated by $\alpha_1$ and $\beta$. Then
$H=\langle\alpha_1,\alpha_1^{-p}\beta^{r+1}\rangle\simeq
C_{p^3}\times C_p$ and $G/H=\langle\alpha_2\rangle\simeq C_{p^2}$.
Let $a=(r+1)^{-1}\in\mathbb Z_{p^2}$, hence
$\beta=\alpha_1^{ap}(\alpha_1^{-p}\beta^{r+1})^a$. Similarly to Case
I, we can define $Y_1,Y_2\in V^*$ such that {\allowdisplaybreaks
\begin{align*}
\alpha_1\ :\ &Y_1\mapsto\zeta_{p^3} Y_1,~Y_2\mapsto Y_2,\\
\alpha_1^{-p}\beta^{r+1}\ :\ &Y_1\mapsto Y_1,~Y_2\mapsto\zeta Y_2,\\
\beta\ :\ &Y_1\mapsto\zeta_{p^2}^a Y_1,~Y_2\mapsto\zeta^a Y_2.
\end{align*}}
Thus $K\cdot Y_1+K\cdot Y_2$ is a representation space of the
subgroup $H$.

Define $x_i=\alpha_2^i\cdot Y_1,y_i=\alpha_2^i\cdot Y_2$ for $0\leq
i\leq p^2-1$. From the relations
$\alpha_1\alpha_2^i=\alpha_2^i\alpha_1\beta^i\beta^{\binom{i}{2}p}$
and $\beta\alpha_2^i=\alpha_2^i\beta^{1+ip}$ it follows that
{\allowdisplaybreaks
\begin{align*}
\alpha_1\ :\ &x_i\mapsto\zeta_{p^3}\zeta_{p^2}^{ai}\zeta^{a\binom{i}{2}} x_i,~ y_i\mapsto\zeta^{ai} y_i\\
\beta\ :\ &x_i\mapsto\zeta_{p^2}^a\zeta^{ai} x_i,~ y_i\mapsto\zeta^a y_i,\\
\alpha_2\ :\ &x_0\mapsto x_1\mapsto\cdots\mapsto x_{p^2-1}\mapsto x_0,\\
&y_0\mapsto y_1\mapsto\cdots\mapsto y_{p^2-1}\mapsto y_0,
\end{align*}}
for $0\leq i\leq p^2-1$.

We find that $Y=(\bigoplus_{0\leq i\leq p^2-1}K\cdot
x_i)\oplus(\bigoplus_{0\leq i\leq p^2-1}K\cdot y_i)$ is a faithful
$G$-subspace of $V^*$. Thus, by Theorem \ref{t2.1}, it suffices to
show that $K(x_i,y_i:0\leq i\leq p^2-1)^{G}$ is rational over $K$.

For $1\leq i\leq p^2-1$, define $u_i=x_i/x_{i-1}$ and
$v_i=y_i/y_{i-1}$. We have now {\allowdisplaybreaks
\begin{align*}
\alpha_1\ :\ &u_i\mapsto\zeta_{p^2}^a\zeta^{a(i-1)} u_i,~ v_i\mapsto\zeta^a v_i,\\
\beta\ :\ &u_i\mapsto\zeta^a u_i,~ v_i\mapsto v_i,\\
\alpha_2\ :\ &u_1\mapsto u_2\mapsto\cdots\mapsto u_{p^2-1}\mapsto (u_1u_2\cdots u_{p^2-1})^{-1},\\
&v_1\mapsto v_2\mapsto\cdots\mapsto v_{p^2-1}\mapsto (v_1v_2\cdots
v_{p^2-1})^{-1},
\end{align*}}
for $1\leq i\leq p^2-1$. From Theorem \ref{t2.2} it follows that if
$K(u_i,v_i:1\leq i\leq p^2-1)^{G}(u,v)$ is rational over $K$, so is
$K(x_i,y_i:0\leq i\leq p^2-1)^{G}$ over $K$.

Since $\beta$ acts in the same way as $\alpha_1^p$ on
$K(u_i,v_i:1\leq i\leq p^2-1)$, we find that $K(u_i,v_i:1\leq i\leq
p^2-1)^{G}=K(u_i,v_i:1\leq i\leq
p^2-1)^{\langle\alpha_1,\alpha_2\rangle}$.

For $1\leq i\leq p^2-1$ define $V_i=v_i/u_i^p$. It follows that
{\allowdisplaybreaks
\begin{align*}
\alpha_1\ :\ &u_i\mapsto\zeta_{p^2}^a\zeta^{a(i-1)} u_i,~ V_i\mapsto V_i,\\
\tag{5.2} \alpha_2\ :\ &u_1\mapsto u_2\mapsto\cdots\mapsto u_{p^2-1}\mapsto (u_1u_2\cdots u_{p^2-1})^{-1},\\
&V_1\mapsto V_2\mapsto\cdots\mapsto V_{p^2-1}\mapsto (V_1V_2\cdots
V_{p^2-1})^{-1},
\end{align*}}
for $1\leq i\leq p^2-1$.

Compare Formula (5.2) with Formula (5.1). They look almost the same.
Apply the proof of Case 1.

\emph{Case III.} $G=\Phi_8(321)c_{p-1}$. Denote by $H$ the abelian
normal subgroup of $G$ generated by $\alpha_1$ and $\beta$. Then
$H\simeq C_{p^2}\times C_{p^2}$ and $G/H\simeq C_{p^2}$. Similarly
to Case I, we can define $Y_1,Y_2\in V^*$ such that
{\allowdisplaybreaks
\begin{align*}
\alpha_1\ :\ &Y_1\mapsto\zeta_{p^2} Y_1,~Y_2\mapsto Y_2,\\
\beta\ :\ &Y_1\mapsto Y_1,~Y_2\mapsto\zeta_{p^2} Y_2.
\end{align*}}
Thus $K\cdot Y_1+K\cdot Y_2$ is a representation space of the
subgroup $H$.

Define $x_i=\alpha_2^i\cdot Y_1,y_i=\alpha_2^i\cdot Y_2$ for $0\leq
i\leq p^2-1$. From the relations
$\alpha_1\alpha_2^i=\alpha_2^i\alpha_1\beta^i\beta^{\binom{i}{2}p}$
and $\beta\alpha_2^i=\alpha_2^i\beta^{1+ip}$ it follows that
{\allowdisplaybreaks
\begin{align*}
\alpha_1\ :\ &x_i\mapsto\zeta_{p^2} x_i,~ y_i\mapsto\zeta_{p^2}^i\zeta^{\binom{i}{2}} y_i\\
\beta\ :\ &x_i\mapsto x_i,~ y_i\mapsto\zeta_{p^2}\zeta^i y_i,\\
\alpha_2\ :\ &x_0\mapsto x_1\mapsto\cdots\mapsto x_{p^2-1}\mapsto x_0,\\
&y_0\mapsto y_1\mapsto\cdots\mapsto y_{p^2-1}\mapsto\zeta y_0,
\end{align*}}
for $0\leq i\leq p^2-1$.

For $1\leq i\leq p^2-1$, define $u_i=x_i/x_{i-1}$ and
$v_i=y_i/y_{i-1}$. We have now {\allowdisplaybreaks
\begin{align*}
\alpha_1\ :\ &u_i\mapsto u_i,~ v_i\mapsto\zeta_{p^2}\zeta^{i-1} v_i,\\
\beta\ :\ &u_i\mapsto u_i,~ v_i\mapsto\zeta v_i,\\
\alpha_2\ :\ &u_1\mapsto u_2\mapsto\cdots\mapsto u_{p^2-1}\mapsto (u_1u_2\cdots u_{p^2-1})^{-1},\\
&v_1\mapsto v_2\mapsto\cdots\mapsto v_{p^2-1}\mapsto\zeta
(v_1v_2\cdots v_{p^2-1})^{-1},
\end{align*}}
for $1\leq i\leq p^2-1$. Since $\beta$ acts in the same way as
$\alpha_1^p$ on $K(u_i,v_i:1\leq i\leq p^2-1)$, we find that
$K(u_i,v_i:1\leq i\leq p^2-1)^{G}=K(u_i,v_i:1\leq i\leq
p^2-1)^{\langle\alpha_1,\alpha_2\rangle}$.

Let $\zeta_{p^3}\in K$ be a primitive $p^3$th root of unity such
that $\zeta_{p^3}^{p^2}=\zeta$. For $1\leq i\leq p^2-1$ define
$w_i=v_i/\zeta_{p^3}$. It follows that {\allowdisplaybreaks
\begin{align*}
\alpha_1\ :\ &u_i\mapsto u_i,~ w_i\mapsto\zeta_{p^2}\zeta^{i-1} w_i,\\
\tag{5.3} \alpha_2\ :\ &u_1\mapsto u_2\mapsto\cdots\mapsto u_{p^2-1}\mapsto (u_1u_2\cdots u_{p^2-1})^{-1},\\
&w_1\mapsto w_2\mapsto\cdots\mapsto w_{p^2-1}\mapsto (w_1w_2\cdots
w_{p^2-1})^{-1},
\end{align*}}
for $1\leq i\leq p^2-1$. Compare Formula (5.3) with Formula (5.1)
(or (5.2)). They look almost the same. Apply the proof of Case 1.

\emph{Case IV.} $G=\Phi_8(222)$. Denote by $H$ the abelian normal
subgroup of $G$ generated by $\alpha_1$ and $\beta$. Then $H\simeq
C_{p^2}\times C_{p^2}$ and $G/H\simeq C_{p^2}$. The proof henceforth
is almost the same as Case III.

\emph{Case V.} $G=\Phi_{14}(321)$. Denote by $H$ the abelian normal
subgroup of $G$ generated by $\alpha_2$ and $\beta$. Then $H\simeq
C_{p^2}\times C_{p^2}$ and $G/H\simeq C_{p^2}$.

As before, we can define $Y_1,Y_2\in V^*$ such that
{\allowdisplaybreaks
\begin{align*}
\alpha_2\ :\ &Y_1\mapsto\zeta_{p^2} Y_1,~Y_2\mapsto Y_2,\\
\beta\ :\ &Y_1\mapsto Y_1,~Y_2\mapsto\zeta_{p^2} Y_2.
\end{align*}}
Thus $K\cdot Y_1+K\cdot Y_2$ is a representation space of the
subgroup $H$.

Define $x_i=\alpha_1^i\cdot Y_1,y_i=\alpha_1^i\cdot Y_2$ for $0\leq
i\leq p^2-1$. From the relations
$\alpha_2\alpha_1^i=\alpha_1^i\alpha_2\beta^{-i}$ it follows that
{\allowdisplaybreaks
\begin{align*}
\alpha_2\ :\ &x_i\mapsto\zeta_{p^2} x_i,~ y_i\mapsto\zeta_{p^2}^{-i} y_i\\
\beta\ :\ &x_i\mapsto x_i,~ y_i\mapsto\zeta_{p^2} y_i,\\
\alpha_1\ :\ &x_0\mapsto x_1\mapsto\cdots\mapsto x_{p^2-1}\mapsto x_0,\\
&y_0\mapsto y_1\mapsto\cdots\mapsto y_{p^2-1}\mapsto\zeta y_0,
\end{align*}}
for $0\leq i\leq p^2-1$.

For $1\leq i\leq p^2-1$, define $u_i=x_i/x_{i-1}$ and
$v_i=y_i/y_{i-1}$. We have now {\allowdisplaybreaks
\begin{align*}
\alpha_2\ :\ &u_i\mapsto u_i,~ v_i\mapsto\zeta_{p^2}^{-1} v_i,\\
\beta\ :\ &u_i\mapsto u_i,~ v_i\mapsto v_i,\\
\alpha_1\ :\ &u_1\mapsto u_2\mapsto\cdots\mapsto u_{p^2-1}\mapsto (u_1u_2\cdots u_{p^2-1})^{-1},\\
&v_1\mapsto v_2\mapsto\cdots\mapsto v_{p^2-1}\mapsto\zeta
(v_1v_2\cdots v_{p^2-1})^{-1},
\end{align*}}
for $1\leq i\leq p^2-1$. Since $\beta$ acts trivially on
$K(u_i,v_i:1\leq i\leq p^2-1)$, we find that $K(u_i,v_i:1\leq i\leq
p^2-1)^{G}=K(u_i,v_i:1\leq i\leq
p^2-1)^{\langle\alpha_1,\alpha_2\rangle}$.

Define $w_1=v_1^{p^2}\zeta^{-1},w_i=v_i/v_{i-1}$ for $2\leq i\leq
p^2-1$. We have now
$K(v_1,\dots,v_{p^2-1})^{\langle\alpha_2\rangle}=K(w_1,\dots,w_{p^2-1})$
and {\allowdisplaybreaks\begin{eqnarray*} \alpha_1&:&w_1\mapsto
w_2^{p^2}w_1,w_2\mapsto w_3\mapsto\cdots\mapsto w_{p^2-1}\mapsto
1/(w_1w_2^{p^2-1}w_3^{p^2-2}\cdots w_{p^2-1}^2).
\end{eqnarray*}}

Define $z_1=w_2,z_i=\alpha_1^{i-1}\cdot w_2$ for $2\leq i\leq
p^2-1$. Then $K(w_i:1\leq i\leq p^2-1)=K(z_i:1\leq i\leq p^2-1)$ and
\begin{align*}
\alpha_1\ :\ &z_1\mapsto z_2\mapsto\cdots\mapsto z_{p^t-1}\mapsto
(z_1z_2\cdots z_{p^2-1})^{-1}.
\end{align*}
The action of $\alpha_1$ can be linearized according to Lemma
\ref{l2.7}. Thus $K(u_i,z_i:1\leq i\leq
p^2-1)^{\langle\alpha_1\rangle}$ is rational over $K$ by Theorem
\ref{t2.1}. We are done.

\emph{Case VI.} $G=\Phi_{14}(222)$. Denote by $H$ the abelian normal
subgroup of $G$ generated by $\alpha_2$ and $\beta$. Then $H\simeq
C_{p^2}\times C_{p^2}$ and $G/H\simeq C_{p^2}$. The proof henceforth
is almost the same as Case V.


\begin{thebibliography}{AAAA}
\bibitem[AHK]{AHK}
H. Ahmad, S. Hajja and M. Kang, Rationality of some projective
linear actions, {\it J. Algebra} {\bf 228} (2000), 643--658.
\bibitem[Be]{Be2}
H. A. Bender, On groups of order $p^m,p$ being an odd prime number,
which contain an abelian subgroup of order $p^{m-1}$, {\it Ann.
Math.}, {\bf 29} No. 1/4 (1927-1928), 88--94.
\bibitem[CK]{CK}
H. Chu and M. Kang, Rationality of $p$-group actions, {\it J.
Algebra} {\bf 237} (2001), 673--690.
\bibitem[GMS]{GMS}
S. Garibaldi, A. Merkurjev and J-P. Serre, ``Cohomological
invariants in Galois cohomology'', AMS Univ. Lecture Series vol. 28,
Amer. Math. Soc., Providence, 2003.
\bibitem[Ha]{Ha}
G.K. Haeuslein, On the invariants of finite groups having an abelian
normal subgroup of prime index, {\it J. London Math. Soc.} {\bf 3}
(1971), 355–-360.
\bibitem[Haj]{Haj}
M. Hajja, Rational invariants of meta-abelian groups of linear
automorphisms, {\it J. Algebra} {\bf 80} (1983), 295–-305.
\bibitem[HK]{HK}
S. Hajja and M. Kang, Some actions of symmetric groups, {\it J.
Algebra} {\bf 177} (1995), 511--535.
\bibitem[HoK]{HoK}
A. Hoshi and M. Kang, Unramified Brauer groups for groups of order
$p^5$, arXiv:1109.2966v1 [math.AC].
\bibitem[HuK]{HuK}
S. J. Hu and M. Kang, Noether's problem for some $p$-groups, in
``Cohomological and geometric approaches to rationality problems'',
edited by F. Bogomolov and Y. Tschinkel, Progress in Math. vol. 282,
Birkh\"auser, Boston, 2010.
\bibitem[Ja]{Ja}
R. James, The groups of order $p^6$ ($p$ an odd prime), {\it Math.
Comp.} {\bf 34} No. 150 (1980), 613--637.
\bibitem[Ka1]{Ka1}
M. Kang, Noether's problem for metacyclic $p$-groups, {\it Adv.
Math.} {\bf 203} (2005), 554--567.
\bibitem[Ka2]{Ka2}
M. Kang, Noether's problem for $p$-groups with a cyclic subgroup of
index $p^2$, {\it Adv. Math.} {\bf 226} (2011) 218--234.
\bibitem[Ka3]{Ka3}
M. Kang, Rationality problem for some meta-abelian groups, {\it J.
Algebra} {\bf 322} (2009), 1214-1219.
\bibitem[Mi]{Mi}
I. Michailov, Noether's problem for $p$-groups with an abelian
subgroup of index $p$, (preprint available at arXiv:1201.5555v3
[math.AG]).
\bibitem[MM]{MM}
J.M. Masley, H.L. Montgomery, Cyclotomic fields with unique
factorization, {\it J. Reine Angew. Math.} {\bf 286/287} (1976)
248-–256.
\bibitem[Sa1]{Sa1}
D. J. Saltman, Generic Galois extensions and problems in field
theory, {\it Adv. Math.} {\bf 43} (1982), 250--283.
\bibitem[Sa2]{Sa2}
D. J. Saltman, Noether's problem over an algebraically closed field,
{\it Invent. Math.} {\bf 77}  (1984), 71--84.
\bibitem[Sw]{Sw}
R. Swan, Noether's problem in Galois theory, in ``Emmy Noether in
Bryn Mawr'', edited by B. Srinivasan and J. Sally, Springer-Verlag,
Berlin, 1983.
\end{thebibliography}
\end{document}